\documentclass[12pt]{article}
\usepackage{latexsym,amssymb,amsmath}
\usepackage{stmaryrd}
\begin{document}
\baselineskip=20pt

\newcommand{\la}{\langle}
\newcommand{\ra}{\rangle}
\newcommand{\psp}{\vspace{0.4cm}}
\newcommand{\pse}{\vspace{0.2cm}}
\newcommand{\ptl}{\partial}
\newcommand{\dlt}{\delta}
\newcommand{\sgm}{\sigma}
\newcommand{\al}{\alpha}
\newcommand{\be}{\beta}
\newcommand{\G}{\Gamma}
\newcommand{\gm}{\gamma}
\newcommand{\vs}{\varsigma}
\newcommand{\Lmd}{\Lambda}
\newcommand{\lmd}{\lambda}
\newcommand{\td}{\tilde}
\newcommand{\vf}{\varphi}
\newcommand{\yt}{Y^{\nu}}
\newcommand{\wt}{\mbox{wt}\:}
\newcommand{\rd}{\mbox{Res}}
\newcommand{\ad}{\mbox{ad}}
\newcommand{\stl}{\stackrel}
\newcommand{\ol}{\overline}
\newcommand{\ul}{\underline}
\newcommand{\es}{\epsilon}
\newcommand{\dmd}{\diamond}
\newcommand{\clt}{\clubsuit}
\newcommand{\vt}{\vartheta}
\newcommand{\ves}{\varepsilon}
\newcommand{\dg}{\dagger}
\newcommand{\tr}{\mbox{Tr}}
\newcommand{\ga}{{\cal G}({\cal A})}
\newcommand{\hga}{\hat{\cal G}({\cal A})}
\newcommand{\Edo}{\mbox{End}\:}
\newcommand{\for}{\mbox{for}}
\newcommand{\kn}{\mbox{ker}}
\newcommand{\Dlt}{\Delta}
\newcommand{\rad}{\mbox{Rad}}
\newcommand{\rta}{\rightarrow}
\newcommand{\mbb}{\mathbb}
\newcommand{\lra}{\Longrightarrow}
\newcommand{\X}{{\cal X}}
\newcommand{\Y}{{\cal Y}}
\newcommand{\Z}{{\cal Z}}
\newcommand{\U}{{\cal U}}
\newcommand{\V}{{\cal V}}
\newcommand{\W}{{\cal W}}
\newcommand{\sta}{\theta}
\setlength{\unitlength}{3pt}

\begin{center}{\Large \bf A New Functor from $D_5$-Mod to $E_6$-Mod} 
\end{center}
\vspace{0.2cm}

\begin{center}{\large Xiaoping Xu
}\end{center}
\begin{center}{Hua Loo-Keng Key Mathematical Laboratory}\end{center}
\begin{center}{Institute of Mathematics, Academy of Mathematics \& System Sciences}\end{center}
\begin{center}{Chinese Academy of Sciences, Beijing 100190, P.R. China
\footnote{Research supported
 by China NSF 11171324}}\end{center}

\begin{abstract}{ We find a new
representation of the simple Lie algebra of type $E_6$ on the
polynomial algebra in 16 variables, which gives a fractional
representation of the corresponding Lie group on 16-dimensional
space. Using this representation and Shen's idea of mixed product,
we construct a functor from $D_5$-{\bf Mod} to $E_6$-{\bf Mod}. A
condition for the functor to map a finite-dimensional irreducible
$D_5$-module to an infinite-dimensional irreducible $E_6$-module is
obtained. Our general frame also gives a direct polynomial extension
from irreducible $D_5$-modules to irreducible $E_6$-modules. The
obtained infinite-dimensional irreducible $E_6$-modules are $({\cal
G},K)$-modules in terms of  Lie group representations. The results
could be used in studying the quantum field theory with $E_6$
symmetry and symmetry of partial differential equations.
}\end{abstract}

\section{Introduction}

\quad \quad A quantum field is an operator-valued function on a
certain Hilbert space, which is often a direct sum of
infinite-dimensional irreducible modules of a certain Lie algebra
(group). The Lie algebra of two-dimensional conformal group is
exactly the Virasoro algebra. The minimal models of two-dimensional
conformal field theory were constructed from direct sums of certain
infinite-dimensional irreducible modules of the Virasoro algebra,
where a distinguished module called, the {\it vacuum module}, gives
rise to a vertex operator algebra.

It is well known that $n$-dimensional projective group gives rise to
a non-homogenous representation of the Lie algebra $sl(n+1,\mbb{C})$
on the polynomial functions of the projective space. Using Shen's
mixed product for Witt algebras, Zhao and the author [ZX]
generalized the above representation of $sl(n+1,\mbb{C})$ to a
non-homogenous representation on the tensor space of any
finite-dimensional irreducible $gl(n,\mbb{C})$-module with the
polynomial space. Moreover, the structure of such a representation
was completely determined by employing projection operator
techniques (cf. [Gm]) and the well-known Kostant's characteristic
identities (cf. [K]). The result can be used to study the quantum
field theory with $sl(n+1,\mbb{C})$ as the symmetry.
 Furthermore, we [XZ] generalize the
conformal representation of of $o(n+2,\mbb{C})$ to a non-homogenous
representation of $o(n+2,\mbb{C})$ on the tensor space of any
finite-dimensional irreducible $o(n,\mbb{C})$-module with a
polynomial space by Shen's idea of mixed product for Witt algebras.
It turns out that  a hidden central transformation is involved. More
importantly,  we find a condition on the constant value taken by the
central transformation  such that the generalized conformal
representation is irreducible. The result would be useful in
higher-dimensional conformal field theory.

This paper is the third work in the program of studying
quantum-field motivated representations of finite-dimensional simple
Lie algebras. It is well known that the minimal dimension of
irreducible modules over the  simple Lie algebra of type $E_6$ is
27. Based on a grading of the simple Lie algebra of type $E_6$, we
find a first-order differential operator representation of the Lie
algebra on the polynomial algebra in 16 independent variables.  In
fact, the corresponding Lie group representation is given by
fractional transformations on 16-dimensional space. Using this
representation and Shen's idea of mixed product, we construct a new
functor from $D_5$-{\bf Mod} to $E_6$-{\bf Mod}, where a hidden
central transformation is involved. More importantly,  a condition
for the functor to map a finite-dimensional irreducible $D_5$-module
to an infinite-dimensional irreducible $E_6$-module in terms of the
constant value taken by the central transformation is obtained. Our
general frame also gives a direct polynomial extension from
irreducible $D_5$-modules to irreducible $E_6$-modules, which can be
applied to obtain explicit bases of irreducible $E_6$-modules from
those of irreducible $D_5$-modules. The result could be useful in
understanding the quantum field theory with $E_6$ symmetry. Our
fractional representation of the $E_6$ Lie group could also be  used
in symmetry analysis of partial differential equations just as the
conformal representation of orthogonal Lie groups do. Our
infinite-dimensional irreducible $E_6$-modules are $({\cal
G},K)$-modules in terms of the corresponding Lie group
representations.

The $E_6$ Lie algebra and group are popular mathematical objects
with broad applications. Dickson [D] (1901) first realized that
there exists an $E_6$-invariant trilinear form on its 27-dimensional
basic irreducible module. The 78-dimensional simple Lie algebra of
type $E_6$ can be realized by all the derivations and multiplication
operators with trace zero on the 27-dimensional exceptional simple
Jordan algebra (e.g., cf. [T], [Ad]). Aschbacher [As] used the
Dickson form to study the subgroup structure of the group $E_6$.
Bion-Nadal [B-N] proved that the $E_6$ Coxeter graph can be realized
as a principal graph of subfactor of the hyperfinite $\Pi_1$ factor.
Brylinski and Kostant [BK] obtained a generalized Capelli identity
on the minimal representation of $E_6$. Binegar and Zierau [BZ]
found a singular representation of $E_6$. Ginzburg [G] proved that
the twisted partial $L$-function on the 27-dimensional
representation of $GE_6(\mbb{C})$ is entire except the points 0 and
1. Iltyakov [I] showed that the field of invariant rational
functions of $E_6$ on the direct sum of finite copies of the basic
module and its dual is purely transcendental. Suzuki and Wakui [SW]
studied the Turaev-Viro-Ocneanu invariant of 3-manifilds derived
from the $E_6$-subfactor. Moreover, Cerchiai and Scotti [CS]
investigated the mapping geometry of the $E_6$ group. Furthermore,
the $(A_2,G_2)$ duality in $E_6$ was obtained by Rubenthaler [R]. In
[X2], the author proved that the space of homogeneous polynomial
solutions with degree $m$ for the dual cubic Dickson invariant
differential operator is exactly a direct sum of $\llbracket m/2
\rrbracket+1$ explicitly determined irreducible $E_6$-submodules and
the whole polynomial algebra is a free module over the polynomial
algebra in the Dickson invariant generated by these solutions.
Moreover, we found in [X3] that the weight matrices of $E_6$ on its
minimal irreducible modules and adjoint modules all generate ternary
orthogonal codes with large minimal distances.

Okamoto and Marshak [OM] constructed a grand unification preson
model with $E_6$ metacolor. The $E_6$ Lie algebra was used in [HH]
to explain the degeneracies encountered in the genetic code as the
result of a sequence of symmetry breakings that have occurred during
its evolution. Wang [W] identified Geoner's model with twisted LG
model and $E_6$ singlets. Morrison, Pieruschka and Wybourne [MPW]
constructed the $E_6$ interacting boson model. Berglund, Candelas et
al. [BCDH] studied instanton contributions to the masses and
couplings of $E_6$ singles. Haba and Matsuoka [HM] found large
lepton flavor mixing in the $E_6$-type unification models.
 Ghezelbash, Shafiekhani and Abolbasani [GSA] derived explicitly a set of Picard-Fuchs
 equations of $N=2$ supersymmetric $E_6$ Yang-Mills
theory. Anderson and Bla$\check z$ek [AB1-AB3] found certain
Clebsch-Gordan coefficients in connection with the $E_6$ unification
model building.  Fern\'{a}ndez-N\'{u}$\td{a}$ez, Garcia-Fuertes and
Perelomov [FGP] used the quantum Calogero-Sutherland model
corresponding to the root system of $E_6$ to calculate
Clebach-Gordan series for this algebra. Howl and King [HK] proposed
a minimal $E_6$ supersymmetric standard model which allows Planck
scale unification, provides a solution to the $\mu$ problem and
predicts a new $Z'$. Das and Laperashvili [DL] studied Preon model
related to family replicated $E_6$ unification.

This work further reveals new beauties of the simple Lie algebra of
type $E_6$. In Section 2, we construct the spin representation of
$o(10,\mbb{C})$ in terms of first-order differential operators on
the polynomial algebra in 16 independent variables from the
lattice-construction of the simple Lie algebra of type $E_6$. We
determine  the decomposition of the polynomial algebra into
irreducible $o(10,\mbb{C})$-submodules in Section 3 by means of
partial differential equations. In Section 4, we realized the simple
Lie algebra of type $E_6$ in terms of first-order differential
operators on the polynomial algebra in 16 independent variables.
Section 5 is devoted to the explicit presentation of the functor
from $D_5$-{\bf Mod} to $E_6$-{\bf Mod}. Finally in Section 6, we
determine a condition  for the functor to map a finite-dimensional
irreducible $D_5$-module to an infinite-dimensional irreducible
$E_6$-module.

 \section{Polynomial Representation of $o(10,\mbb{C})$ via $E_6$}

We start with the root lattice construction of the simple Lie
algebra of type $E_6$. As we all know, the Dynkin diagram of $E_6$
is as follows:

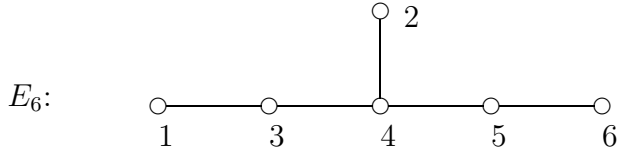
\begin{picture}(80,23)
\put(2,0){$E_6$:}\put(21,0){\circle{2}}\put(21,
-5){1}\put(22,0){\line(1,0){12}}\put(35,0){\circle{2}}\put(35,
-5){3}\put(36,0){\line(1,0){12}}\put(49,0){\circle{2}}\put(49,
-5){4}\put(49,1){\line(0,1){10}}\put(49,12){\circle{2}}\put(52,10){2}\put(50,0){\line(1,0){12}}
\put(63,0){\circle{2}}\put(63,-5){5}\put(64,0){\line(1,0){12}}\put(77,0){\circle{2}}\put(77,
-5){6}
\end{picture}
\vspace{0.7cm}

\noindent For convenience, we will use the notion
$\ol{i,i+j}=\{i,i+1,i+2,...,i+j\}$ for integer $i$ and positive
integer $j$ throughout this paper. Let $\{\al_i\mid i\in\ol{1,6}\}$
be the simple positive roots corresponding to the vertices in the
diagram, and let $\Phi_{E_6}$ be the root system of $E_6$. Set
$$Q_{E_6}=\sum_{i=1}^6\mbb{Z}\al_i,\eqno(2.1)$$ the root lattice of type
$E_6$. Denote by $(\cdot,\cdot)$ the symmetric $\mbb{Z}$-bilinear
form on $Q_{E_6}$ such that
$$\Phi_{E_6}=\{\al\in Q_{E_6}\mid (\al,\al)=2\}.\eqno(2.2)$$
 Define a map $F:
Q_{E_6}\times Q_{E_6}\rta \{1,-1\}$ by
$$F(\sum_{i=1}^6k_i\al_i,\sum_{j=1}^6l_j\al_j)=(-1)^{\sum_{i=1}^6k_il_i+
k_1l_3+k_4l_2+k_3l_4+k_5l_4+k_6l_5},\qquad
k_i,l_j\in\mbb{Z}.\eqno(2.3)$$ Then for $\al,\be,\gm\in Q_{E_6}$,
$$F(\al+\be,\gm)=F(\al,\gm)F(\be,\gm),\;\;F(\al,\be+\gm)=F(\al,\be)F(\al,\gm),
\eqno(2.4)$$
$$F(\al,\be)F(\be,\al)^{-1}=(-1)^{(\al,\be)},\;\;F(\al,\al)=(-1)^{(\al,\al)/2}.
\eqno(2.5)$$ In particular,
$$F(\al,\be)=-F(\be,\al)\qquad
\mbox{if}\;\;\al,\be,\al+\be\in\Phi_{E_6}.\eqno(2.6)$$

Denote
$$H=\bigoplus_{i=1}^6\mbb{C}\al_i\eqno(2.7)$$
 and $\mbb{C}$-bilinearly extend $(\cdot,\cdot)$ on $H$. Then the simple  Lie algebra of type $E_6$ is
$${\cal
G}^{E_6}=H\oplus\bigoplus_{\al\in\Phi_{E_6}}\mbb{C}E_{\al}\eqno(2.8)$$
 with the Lie bracket $[\cdot,\cdot]$ determined by:
 $$[H,H]=0,\;\;[h,E_{\al}]=-[E_{\al},h]=(h,\al)E_{\al},\;\;[E_{\al},E_{-\al}]=-\al,
 \eqno(2.9)$$
 $$[E_{\al},E_{\be}]=\left\{\begin{array}{ll}0&\mbox{if}\;\al+\be\not\in\Phi_{E_6},\\
 F(\al,\be)E_{\al+\be}&\mbox{if}\;\al+\be\in\Phi_{E_6}\end{array}\right.\eqno(2.10)$$(e.g., cf. [Ka],
[X1]). Moreover, we define a bilinear form $(\cdot|\cdot)$ on ${\cal
G}^{E_6}$ by
$$(h_1|h_2)=(h_1,h_2),\;\; (h|E_{\al})=0,\;\;
 (E_{\al}|E_{\be})=-\dlt_{\al+\be,0}\eqno(2.11)$$
for $h_1,h_2\in H$ and $\al,\be\in \Phi_{E_6}$. It can be verified
that $(\cdot|\cdot)$ is a ${\cal G}^{E_6}$-invariant form, that is,
$$([u,v]|w)=(u|[v,w])\qquad\for\;\;u,v\in{\cal
G}^{E_6}.\eqno(2.12)$$

Let
$$Q^{D_5}=\sum_{i=1}^5\mbb{Z}\al_i,\qquad \Phi_{D_5}=\Phi_{E_6}\bigcap Q^{D_5}. \eqno(2.13)$$
Then
$${\cal G}^{D_5}=\sum_{i=1}^5\mbb{C}\al_i+\sum_{\be\in
\Phi_{D_5}}\mbb{C}E_\be\eqno(2.14)$$ forms a Lie subalgebra of
${\cal G}^{E_6}$, which is isomorphic to the orthogonal Lie algebra
\begin{eqnarray*}\qquad o(10,\mbb{C})&=&\sum_{1\leq
p<q\leq 5}[\mbb{C}(E_{p,n+q}-E_{q,n+p})+
\mbb{C}(E_{n+p,q}-E_{n+q,p})]\\ &
&+\sum_{i,j=1}^5\mbb{C}(E_{i,j}-E_{n+j,n+i}).\hspace{7.1cm}(2.15)\end{eqnarray*}
 Denote by
$\Phi_{E_6}^+$ the set of positive roots of $E_6$ and by
$\Phi_{D_5}^+$ the set of positive roots of $D_5$. We find the
elements of $\Phi_{D_5}^+$:
$$\al_r\;(r\in\ol{1,5}),\;\al_1+\al_3,\;\al_2+\al_4,\;\al_3+\al_4,\;\al_4+\al_5,\;
\al_1+\al_3+\al_4,\;\al_2+\al_3+\al_4,\eqno(2.16)$$
$$\al_2+\al_4+\al_5,\;\al_3+\al_4+\al_5,\sum_{r=1}^4\al_r,\;\al_1+\al_3+\al_4+\al_5,\eqno(2.17)$$
$$\sum_{i=2}^5\al_i,\;
\sum_{s=1}^5\al_s,\;\al_4+\sum_{i=2}^5\al_i,\;\al_4+\sum_{i=1}^5\al_i,\;
\al_3+\al_4+\sum_{i=1}^5\al_i.\eqno(2.18)$$ Moreover, the elements
in $\Phi_{E_6}^+\setminus\Phi_{D_5}^+$ are:
$$\al_6,\;\al_5+\al_6,\;\sum_{r=4}^6\al_r,\;\sum_{i=3}^6\al_i,\;\al_2+\sum_{r=4}^6\al_r,\;\sum_{i=2}^6\al_i,\;
\al_1+\sum_{i=3}^6\al_i,\;\sum_{i=1}^6\al_i,\;\al_4+\sum_{i=2}^6\al_i,\eqno(2.19)$$
$$\al_4+\sum_{i=1}^6\al_i,\;\al_4+\al_5+\sum_{i=2}^6\al_i,\;\al_3+\al_4+\sum_{i=1}^6\al_i,\;
\al_4+\al_5+\sum_{i=1}^6\al_i,\eqno(2.20)$$
$$\sum_{i=1}^6\al_i+\sum_{r=3}^5\al_r,\;
\al_4+\sum_{i=1}^6\al_i+\sum_{r=3}^5\al_r,\;\al_4+\sum_{i=1}^6\al_i+\sum_{r=2}^5\al_r.
\eqno(2.21)$$

For convenience, we denote
$$\xi_1=E_{\al_6},\;\;\xi_2=E_{\al_5+\al_6}
,\;\;\xi_3=E_{\sum_{r=4}^6\al_r},\;\;\xi_4=E_{\sum_{i=3}^6\al_i},
\eqno(2.22)$$
$$\xi_5=E_{\al_2+\sum_{r=4}^6\al_r},\;\;
\xi_6=E_{\sum_{i=2}^6\al_i},\;\;
\xi_7=E_{\al_1+\sum_{i=3}^6\al_i},\;\;\xi_8=E_{\sum_{i=1}^6\al_i},\eqno(2.23)$$
$$\xi_9=E_{\al_4+\sum_{i=2}^6\al_i},\;\;\xi_{10}=E_{\al_4+\sum_{i=1}^6\al_i},\;
\xi_{11}=E_{\al_4+\al_5+\sum_{i=2}^6\al_i},\eqno(2.24)$$
$$ \xi_{12}=E_{\al_3+\al_4
+\sum_{i=1}^6\al_i},\;\;\xi_{13}=E_{\al_4+\al_5+\sum_{i=1}^6\al_i},\;\;
\xi_{14}=E_{\sum_{i=1}^6\al_i+\sum_{r=3}^5\al_r},\eqno(2.25)$$
$$\xi_{15}=E_{\al_4+\sum_{i=1}^6\al_i+\sum_{r=3}^5\al_r},\;\;
\xi_{16}=E_{\al_4+\sum_{i=1}^6\al_i+\sum_{r=2}^5\al_r},
\eqno(2.26)$$
$$\eta_1=E_{-\al_6},\;\;\eta_2=E_{-\al_5-\al_6},\;\;\eta_3=E_{-\sum_{r=4}^6\al_r},\;\;\eta_4=E_{-\sum_{i=3}^6\al_i},
\eqno(2.27)$$
$$\eta_5=E_{-\al_2-\sum_{r=4}^6\al_r},\;\;
\eta_6=E_{-\sum_{i=2}^6\al_i},\;\;
\eta_7=E_{-\al_1-\sum_{i=3}^6\al_i},\;\;\eta_8=E_{-\sum_{i=1}^6\al_i},\eqno(2.28)$$
$$\eta_9=E_{-\al_4-\sum_{i=2}^6\al_i},\;\;\eta_{10}=E_{-\al_4-\sum_{i=1}^6\al_i},\;
\eta_{11}=E_{-\al_4-\al_5-\sum_{i=2}^6\al_i},\eqno(2.29)$$
$$ \eta_{12}=E_{-\al_3-\al_4
-\sum_{i=1}^6\al_i},\;\;\eta_{13}=E_{-\al_4-\al_5-\sum_{i=1}^6\al_i},\;\;
\eta_{14}=E_{-\sum_{i=1}^6\al_i-\sum_{r=3}^5\al_r},\eqno(2.30)$$
$$\eta_{15}=E_{-\al_4-\sum_{i=1}^6\al_i-\sum_{r=3}^5\al_r},\;\;\eta_{16}=E_{-\al_4-\sum_{i=1}^6\al_i-\sum_{r=2}^5\al_r}.
\eqno(2.31)$$ Set
$${\cal G}_+=\sum_{i=1}^{16}\mbb{C}\xi_i,\qquad{\cal
G}_-=\sum_{i=1}^{16}\mbb{C}\eta_i,\qquad{\cal G}_0={\cal
G}^{D_5}+\mbb{F}\al_6.\eqno(2.32)$$ It is straightforward to verify
that ${\cal G}_\pm$ are ableian Lie subalgebras of ${\cal G}^{E_6}$,
${\cal G}_0$ is a reductive Lie subalgebra of ${\cal G}^{E_6}$ and
$${\cal G}^{E_6}={\cal G}_-\oplus{\cal G}_0\oplus{\cal
G}_+.\eqno(2.33)$$ Moreover, ${\cal G}_\pm$ form irreducible ${\cal
G}_0$-submodules with respect to the adjoint representation of
${\cal G}^{E_6}$. Furthermore,
$$(\xi_i|\eta_j)=-\dlt_{i,j}\qquad\for\;\;i,j\in\ol{1,16}\eqno(2.34)$$
by (2.11). Expression (2.12) shows that ${\cal G}_+$ is isomorphic
to the dual ${\cal G}_0$-module of ${\cal G}_-$.

Set
$${\cal A}=\mbb{C}[x_1,x_2,...,x_{16}],\eqno(2.35)$$
the polynomial algebra in $x_1,x_2,...,x_{16}$. Write
$$[u,\eta_i]=\sum_{j=1}^{16}\vf_{i,j}(u)\eta_j\qquad\for\;\;i\in\ol{1,16},\;u\in{\cal
G}_0,\eqno(2.36)$$ where $\vf_{i,j}(u)\in\mbb{C}$. Define an action
of ${\cal G}_0$ on ${\cal A}$ by
$$u(f)=\sum_{i,j=1}^{16}\vf_{i,j}(u)x_j\ptl_{x_i}(f)\qquad\for\;\;u\in{\cal
G}_0,\;f\in{\cal A}.\eqno(2.37)$$ Then ${\cal A}$ forms a ${\cal
G}_0$-module and the subspace
$$V=\sum_{i=1}^{16}\mbb{C}x_i\eqno(2.38)$$
forms a ${\cal G}_0$-submodule isomorphic to ${\cal G}_-$, where the
isomorphism is determined by $x_i\mapsto \eta_i$  for
$i\in\ol{1,16}$.

Denote by $\mbb{N}$ the set of nonnegative integers. Write
$$x^\al=\prod_{i=1}^{16}x_i^{\al_i},\;\;\ptl^\al=\prod_{i=1}^{16}\ptl_{x_i}^{\al_i}\qquad\for\;\;
\al=(\al_1,\al_2,...,\al_{16})\in\mbb{N}^{16}.\eqno(2.39)$$ Let
$$\mbb{A}=\sum_{\al\in\mbb{N}^{16}}{\cal A}\ptl^\al\eqno(2.40)$$
be the algebra of differential operators on ${\cal A}$. Then the
linear transformation $\tau$ determined by
$$\tau(x^\be\ptl^\gm)=x^\gm\ptl^\be\qquad\for\;\;\be,\gm\in\mbb{N}^{16}\eqno(2.41)$$
is an involutive anti-automorphism of $\mbb{A}$.

 According
to (2.9) and (2.10), we have the Lie algebra isomorphism $\nu:
o(10,\mbb{C})\rta {\cal G}^{D_5}$ determined by the generators:
$$\nu(E_{1,2}-E_{7,6})=E_{\al_1},\;\;\nu(E_{2,3}-E_{8,7})=E_{\al_3},\;\;
\nu(E_{3,4}-E_{9,8})=E_{\al_4},\eqno(2.42)$$
$$\nu(E_{4,5}-E_{10,9})=E_{\al_5},\;\;\nu(E_{4,10}-E_{5,9})=E_{\al_2},\;\;
\nu(E_{2,1}-E_{6,7})=-E_{-\al_1}\eqno(2.43)$$
$$\nu(E_{3,2}-E_{7,8})=-E_{-\al_3},\;\;\nu(E_{4,3}-E_{8,9})=-E_{-\al_4},\;\;
\nu(E_{5,4}-E_{9,10})=-E_{-\al_5},\eqno(2.44)$$
$$\nu(E_{10,4}-E_{9,5})=-E_{-\al_2},\;\;\nu(E_{1,1}-E_{6,6})=\al_1+\al_3+\al_4+\frac{1}{2}(
\al_2+\al_5),\eqno(2.45)$$
$$\nu(E_{2,2}-E_{7,7})=\al_3+\al_4+\frac{1}{2}(
\al_2+\al_5),\;\;\nu(E_{3,3}-E_{8,8})=\al_4+\frac{1}{2}(
\al_2+\al_5),\eqno(2.46)$$
$$\nu(E_{4,4}-E_{9,9})=\frac{1}{2}(
\al_2+\al_5),\qquad\nu(E_{5,5}-E_{10,10})=\frac{1}{2}(
\al_2-\al_5).\eqno(2.47)$$ Then ${\cal A}$ becomes an
$o(10,\mbb{C})$-module with respect to the action
$$A(f)=\nu(A)(f)\qquad\for\;\;A\in o(10,\mbb{C}),\;f\in\mbb{\cal A}
.\eqno(2.48)$$

Thanks to (2.9), (2.10), (2.27)-(2.31), (2.36) and (2.37), we have
$$(E_{1,2}-E_{7,6})|_{\cal A}=x_4\ptl_{x_7}+x_6\ptl_{x_8}+x_9\ptl_{x_{10}}+x_{11}\ptl_{x_{13}},\eqno(2.49)$$
$$(E_{2,3}-E_{8,7})|_{\cal A}=x_3\ptl_{x_4}+x_5\ptl_{x_6}+x_{10}\ptl_{x_{12}}+x_{13}\ptl_{x_{14}},\eqno(2.50)$$
$$(E_{3,4}-E_{9,8})|_{\cal A}=-x_2\ptl_{x_3}-x_6\ptl_{x_9}-x_8\ptl_{x_{10}}+x_{14}\ptl_{x_{15}},\eqno(2.51)$$
$$(E_{4,5}-E_{10,9})|_{\cal A}=-x_1\ptl_{x_2}+x_9\ptl_{x_{11}}+x_{10}\ptl_{x_{13}}+x_{12}\ptl_{x_{14}},\eqno(2.52)$$
$$(E_{4,10}-E_{5,9})|_{\cal A}=-x_3\ptl_{x_5}-x_4\ptl_{x_6}-x_7\ptl_{x_8}+x_{15}\ptl_{x_{16}},\eqno(2.53)$$
$$(E_{1,3}-E_{8,6})|_{\cal A}=-x_3\ptl_{x_7}-x_5\ptl_{x_8}+x_9\ptl_{x_{12}}+x_{11}\ptl_{x_{14}},\eqno(2.54)$$
$$(E_{2,4}-E_{9,7})|_{\cal A}=x_2\ptl_{x_4}-x_5\ptl_{x_9}+x_8\ptl_{x_{12}}+x_{13}\ptl_{x_{15}},\eqno(2.55)$$
$$(E_{3,5}-E_{10,8})|_{\cal A}=-x_1\ptl_{x_3}-x_6\ptl_{x_{11}}-x_8\ptl_{x_{13}}-x_{12}\ptl_{x_{15}},\eqno(2.56)$$
$$(E_{3,10}-E_{5,8})|_{\cal A}=x_2\ptl_{x_5}-x_4\ptl_{x_9}-x_7\ptl_{x_{10}}+x_{14}\ptl_{x_{16}},\eqno(2.57)$$
$$(E_{1,4}-E_{9,6})|_{\cal A}=-x_2\ptl_{x_7}+x_5\ptl_{x_{10}}+x_6\ptl_{x_{12}}+x_{11}\ptl_{x_{15}},\eqno(2.58)$$
$$(E_{2,5}-E_{10,7})|_{\cal A}=x_1\ptl_{x_4}-x_5\ptl_{x_{11}}+x_8\ptl_{x_{14}}-x_{10}\ptl_{x_{15}},\eqno(2.59)$$
$$(E_{2,10}-E_{5,7})|_{\cal A}=-x_2\ptl_{x_6}-x_3\ptl_{x_9}+x_7\ptl_{x_{12}}+x_{13}\ptl_{x_{16}},\eqno(2.60)$$
$$(E_{3,9}-E_{4,8})|_{\cal A}=-x_1\ptl_{x_5}+x_4\ptl_{x_{11}}+x_7\ptl_{x_{13}}+x_{12}\ptl_{x_{16}},\eqno(2.61)$$
$$(E_{1,5}-E_{10,6})|_{\cal A}=-x_1\ptl_{x_7}+x_5\ptl_{x_{13}}+x_6\ptl_{x_{14}}-x_9\ptl_{x_{15}},\eqno(2.62)$$
$$(E_{1,10}-E_{5,6})|_{\cal A}=x_2\ptl_{x_8}+x_3\ptl_{x_{10}}+x_4\ptl_{x_{12}}+x_{11}\ptl_{x_{16}},\eqno(2.63)$$
$$(E_{2,9}-E_{4,7})|_{\cal A}=x_1\ptl_{x_6}+x_3\ptl_{x_{11}}-x_7\ptl_{x_{14}}+x_{10}\ptl_{x_{16}},\eqno(2.64)$$
$$(E_{1,9}-E_{4,6})|_{\cal A}=-x_1\ptl_{x_8}-x_4\ptl_{x_{14}}-x_3\ptl_{x_{13}}+x_9\ptl_{x_{16}},\eqno(2.65)$$
$$(E_{2,8}-E_{3,7})|_{\cal A}=x_1\ptl_{x_9}-x_2\ptl_{x_{11}}+x_7\ptl_{x_{15}}-x_8\ptl_{x_{16}},\eqno(2.66)$$
$$(E_{1,8}-E_{3,6})|_{\cal A}=-x_1\ptl_{x_{10}}+x_2\ptl_{x_{13}}+x_4\ptl_{x_{15}}-x_6\ptl_{x_{16}},\eqno(2.67)$$
$$(E_{1,7}-E_{2,6})|_{\cal A}=x_1\ptl_{x_{12}}-x_2\ptl_{x_{14}}+x_3\ptl_{x_{15}}-x_5\ptl_{x_{16}},\eqno(2.68)$$
$$(E_{j,i}-E_{5+i,5+j})|_{\cal A}=\tau[(E_{i,j}-E_{5+j,5+i})|_{\cal A}],\eqno(2.69)$$
$$(E_{5+j,i}-E_{5+i,j})|_{\cal A}=\tau[(E_{i,5+j}-E_{j,5+i})|_{\cal A}]\eqno(2.70)$$
for $1\leq i<j\leq 16$,
$$(E_{r,r}-E_{5+r,5+r})|_{\cal A}=\sum_{i=1}^{16}(1/2+a_{r,i})x_i\ptl_{x_i},\qquad
r\in\ol{1,5},\eqno(2.71)$$ where $a_{r,i}$ are given in the
following table
\begin{center}{\bf \large Table 1}\end{center}
\begin{center}\begin{tabular}{|r||r|r|r|r|r|r|r|r|r|r|r|r|r|r|r|r|}\hline
$i$&1&2&3&4&5&6&7&8&9&10&11&12&13&14&15&16
\\\hline\hline
$a_{1,i}$&0&0&0&0&0&0&$-1$&$-1$&0&$-1$&0&$-1$&$-1$&$-1$&$-1$&$-1$
\\\hline
$a_{2,i}$&0&0&0&$-1$&0&$-1$&0&0&$-1$&0&$-1$&$-1$&0&$-1$&$-1$&$-1$
\\\hline$a_{3,i}$&0&0&$-1$&0&$-1$&0&0&0&$-1$&$-1$&$-1$&0&$-1$&0&$-1$&$-1$
\\\hline$a_{4,i}$&0&$-1$&0&0&$-1$&$-1$&0&$-1$&0&0&$-1$&0&$-1$&$-1$&0&$-1$
\\\hline$a_{5,i}$&$-1$&0&0&0&$-1$&$-1$&0&$-1$&$-1$&$-1$&0&$-1$&0&0&0&$-1$\\\hline
\end{tabular}\end{center}
Note that (2.39)-(2.68) are the representation formulas of all the
positive root vectors. In particular, $x_1$ is a highest-weight
vector of $V$ with weight $\lmd_4$, the forth fundamental weight of
$E_6$, and $V$ gives a spin representation of $o(10,\mbb{C})$.

\section{Decomposition of the $o(10,\mbb{C})$-Module ${\cal A}$}

Recall that the representation of ${\cal G}^{D_5}$ on ${\cal A}$ is
given by (2.37) and the representation of $o(10,\mbb{C})$ is given
in (2.48). We calculate
$$\al_r|_{\cal
A}=\sum_{i=1}^{16}b_{r,i}x_i\ptl_{x_i}\qquad\for\;\;r\in\ol{1,5},\eqno(3.1)$$
where $b_{r,i}$ are given in the following table:
\begin{center}{\bf \large Table 2}\end{center}
\begin{center}\begin{tabular}{|r||r|r|r|r|r|r|r|r|r|r|r|r|r|r|r|r|}\hline
$i$&1&2&3&4&5&6&7&8&9&10&11&12&13&14&15&16
\\\hline\hline
$b_{1,i}$&0&0&0&1&0&1&$-1$&$-1$&1&$-1$&1&0&$-1$&0&0&0
\\\hline
$b_{2,i}$&0&0&1&1&$-1$&$-1$&1&$-1$&0&0&0&0&0&0&1&$-1$
\\\hline$b_{3,i}$&0&0&1&$-1$&1&$-1$&0&0&0&1&0&$-1$&1&$-1$&0&0
\\\hline$b_{4,i}$&0&1&$-1$&0&0&1&0&1&$-1$&$-1$&0&0&0&1&$-1$&0
\\\hline$b_{5,i}$&1&$-1$&0&0&0&0&0&0&1&1&$-1$&1&$-1$&$-1$&0&0\\\hline
\end{tabular}\end{center}
Recall that  a singular vector of $o(10,\mbb{C})$ is a nonzero
weight vector annihilated by positive root vectors. Note that the
weight of a singular vector in ${\cal A}$ must be dominate integral.
The above table motivates us to assume that
$$\zeta_1=a_1x_1x_{11}+a_2x_2x_9+a_3x_3x_6+a_4x_4x_5\eqno(3.2)$$
is a singular vector, where $a_i$ are constants to be determined. By
(2.49),
$$(E_{1,2}-E_{7,6})(\zeta_1)=0.\eqno(3.3)$$
Moreover, (2.50) says
$$(E_{2,3}-E_{8,7})(\zeta_1)=(a_3+a_4)x_3x_5=0\lra
a_4=-a_3.\eqno(3.4)$$ Expression (2.51) gives
$$(E_{3,4}-E_{9,8})(\zeta_1)=-(a_2+a_3)x_2x_6=0\lra
a_3=-a_2.\eqno(3.5)$$ Furthermore, (2.52) yields
$$(E_{4,5}-E_{10,9})(\zeta_1)=(a_1-a_2)x_1x_9=0\lra
a_2=a_1.\eqno(3.6)$$ According to (2.53) and (3.4),
$$(E_{4,10}-E_{5,9})(\zeta_1)=-(a_3+a_4)x_3x_4=0.\eqno(3.7)$$
Taking $a_1=1$, we have the singular vector
$$\zeta_1=x_1x_{11}+x_2x_9-x_3x_6+x_4x_5\eqno(3.8)$$
of weight $\lmd_1$, the first fundamental weight of $E_6$. Thus
$\zeta_1$ generates the 10-dimensional natural
$o(10,\mbb{C})$-module $U$. According to (2.49)-(2.53), (2.69) and
(2.70),
$$(E_{2,1}-E_{6,7})|_{\cal A}=x_7\ptl_{x_4}+x_8\ptl_{x_6}+x_{10}\ptl_{x_9}+x_{13}\ptl_{x_{11}},\eqno(3.9)$$
$$(E_{3,2}-E_{7,8})|_{\cal A}=x_4\ptl_{x_3}+x_6\ptl_{x_5}+x_{12}\ptl_{x_{10}}+x_{14}\ptl_{x_{13}},\eqno(3.10)$$
$$(E_{4,3}-E_{8,9})|_{\cal A}=-x_3\ptl_{x_2}-x_9\ptl_{x_6}-x_{10}\ptl_{x_8}+x_{15}\ptl_{x_{14}},\eqno(3.11)$$
$$(E_{5,4}-E_{9,10})|_{\cal A}=-x_2\ptl_{x_1}+x_{11}\ptl_{x_9}+x_{13}\ptl_{x_{10}}+x_{14}\ptl_{x_{12}},\eqno(3.12)$$
$$(E_{10,4}-E_{9,5})|_{\cal A}=-x_5\ptl_{x_3}-x_6\ptl_{x_4}-x_8\ptl_{x_7}+x_{16}\ptl_{x_{15}}.\eqno(3.13)$$
We take
$$\zeta_2=(E_{2,1}-E_{6,7})(\zeta_1)=x_1x_{13}+x_2x_{10}-x_3x_8+x_5x_7,\eqno(3.14)$$
$$\zeta_3=(E_{3,2}-E_{7,8})(\zeta_2)=x_1x_{14}+x_2x_{12}-x_4x_8+x_6x_7,\eqno(3.15)$$
$$\zeta_4=(E_{4,3}-E_{8,9})(\zeta_3)=x_1x_{15}-x_3x_{12}+x_4x_{10}-x_7x_9,\eqno(3.16)$$
$$\zeta_5=(E_{5,4}-E_{9,10})(\zeta_4)=-x_2x_{15}-x_3x_{14}+x_4x_{13}-x_7x_{11},\eqno(3.17)$$
$$\zeta_{10}=(E_{10,4}-E_{9,5})(\zeta_4)=x_1x_{16}+x_5x_{12}-x_6x_{10}+x_8x_9,\eqno(3.18)$$
$$\zeta_9=(E_{9,5}-E_{10,4})(\zeta_5)=x_2x_{16}-x_5x_{14}+x_6x_{13}-x_8x_{11},\eqno(3.19)$$
$$\zeta_8=(E_{8,9}-E_{4,3})(\zeta_9)=x_3x_{16}+x_5x_{15}+x_9x_{13}-x_{10}x_{11},\eqno(3.20)$$
$$\zeta_7=(E_{7,8}-E_{3,2})(\zeta_8)=-x_4x_{16}-x_6x_{15}-x_9x_{14}+x_{11}x_{12},\eqno(3.21)$$
$$\zeta_6=(E_{6,7}-E_{2,1})(\zeta_7)=x_7x_{16}+x_8x_{15}+x_{10}x_{14}-x_{12}x_{13}.\eqno(3.22)$$
Then $U=\sum_{i=1}^{10}\mbb{C}\zeta_i$ forms an
$o(10,\mbb{C})$-module isomorphic to the 10-dimensional natural
$o(10,\mbb{C})$-module with $\{\zeta_1,...,\zeta_{10}\}$ as the
standard basis.\psp

{\bf Theorem 3.1}. {\it Any singular vector is a polynomial in $x_1$
and $\zeta_1$.}

{\it Proof}. Note
$$x_{11}=x_1^{-1}(\zeta_1-x_2x_9+x_3x_6-x_4x_5),\eqno(3.23)$$
$$x_{13}=x_1^{-1}(\zeta_2-x_2x_{10}+x_3x_8-x_5x_7),\eqno(3.24)$$
$$x_{14}=x_1^{-1}(\zeta_3-x_2x_{12}+x_4x_8-x_6x_7),\eqno(3.25)$$
$$x_{15}=x_1^{-1}(\zeta_4+x_3x_{12}-x_4x_{10}+x_7x_9),\eqno(3.26)$$
$$x_{16}=x_1^{-1}(\zeta_{10}-x_5x_{12}+x_6x_{10}-x_8x_9).\eqno(3.27)$$
Let $f$ be a singular vector in ${\cal A}$. Substituting
(3.23)-(3.27) into it, we can write
$$f=g(x_i,\zeta_1,\zeta_2,\zeta_3,\zeta_4,\zeta_{10}\mid 11,13,14,15,16\neq
i\ol{1,16}).\eqno(3.28)$$ By (2.52), (2.56), (2.59), (2.61), (2.62)
and (2.64)-(2.68),
$$(E_{4,5}-E_{10,9})(f)=-x_1\ptl_{x_2}(g)=0\lra g_{x_2}=0,\eqno(3.29)$$
$$(E_{3,5}-E_{10,8})(f)=-x_1\ptl_{x_3}(g)=0\lra g_{x_3}=0,\eqno(3.30)$$
$$(E_{2,5}-E_{10,7})(f)=x_1\ptl_{x_4}(g)=0\lra g_{x_4}=0,\eqno(3.31)$$
$$(E_{3,9}-E_{4,8})(f)=-x_1\ptl_{x_5}(g)=0\lra g_{x_5}=0,\eqno(3.32)$$
$$(E_{1,5}-E_{10,6})(f)=-x_1\ptl_{x_7}(g)=0\lra g_{x_7}=0,\eqno(3.33)$$
$$(E_{2,9}-E_{4,7})(f)=x_1\ptl_{x_6}(g)=0\lra g_{x_6}=0,\eqno(3.34)$$
$$(E_{1,9}-E_{4,6})(f)=-x_1\ptl_{x_8}(g)=0\lra g_{x_8}=0,\eqno(3.35)$$
$$(E_{2,8}-E_{3,7})(f)=x_1\ptl_{x_9}(g)=0\lra g_{x_9}=0,\eqno(3.36)$$
$$(E_{1,8}-E_{3,6})(f)=-x_1\ptl_{x_{10}}(g)=0\lra g_{x_{10}}=0,\eqno(3.37)$$
$$(E_{1,7}-E_{2,6})(f)=x_1\ptl_{x_{12}}(g)=0\lra g_{x_{12}}=0.\eqno(3.38)$$
Thus $f$ is a function in $x_1,\zeta_1,\zeta_2,\zeta_3,\zeta_4$ and
$\zeta_{10}$.

According to (2.49)-(2.51) and (2.53),
$$(E_{1,2}-E_{7,6})(f)=\zeta_1\ptl_{\zeta_2}(g)=0\lra g_{\zeta_2}=0,\eqno(3.39)$$
$$(E_{2,3}-E_{8,7})(f)=\zeta_2\ptl_{\zeta_3}(g)=0\lra g_{\zeta_3}=0,\eqno(3.40)$$
$$(E_{3,4}-E_{9,8})(f)=\zeta_3\ptl_{\zeta_4}(g)=0\lra g_{\zeta_4}=0,\eqno(3.41)$$
$$(E_{4,10}-E_{5,9})=\zeta_4\ptl_{\zeta_{10}}(g)=0\lra g_{\zeta_{10}}=0.\eqno(3.42)$$
Hence $f$ is a function in $x_1$ and $\zeta_1$. Thanks to (3.23), it
must be a polynomial in $x_1$ and $\zeta_1.\qquad\Box$\psp

Let $V_{m_1,m_2}$ be the irreducible $o(10,\mbb{C})$-submodule
generated by $x_1^{m_1}\zeta^{m_2}$. By Weyl's theorem of complete
reducibility,
$${\cal A}=\bigoplus_{m_1,m_2=0}^\infty V_{m_1,m_2}.\eqno(3.43)$$
Denote by $V(\lmd)$ the highest-weight irreducible
$o(10,\mbb{C})$-module with the highest weight $\lmd$. The above
equation leads to the following combinatorial identity:
$$\sum_{m_1,m_2=0}^\infty (\dim
V(m_2\lmd_1+m_1\lmd_4))q^{m_1+2m_2}=\frac{1}{(1-q)^{16}},\eqno(3.44)$$
which was proved in (3.23)-(3.41)  by partial differential
equations.

\section{Realization of $E_6$ in 16-Dimensional Space}

 In this section, we
want to find a differential-operator representation of ${\cal
G}^{E_6}$, equivalently, a fraction representation on 16-dimensional
space of the Lie group of type $E_6$.

According to (2.36) and (2.37), we calculate
$$\al_6|_{\cal A}=-2x_1\ptl_{x_1}-\sum_{i=2}^{10}x_i\ptl_{x_i}-x_{12}\ptl_{x_{12}}.\eqno(4.1)$$
 Write
$$\widehat\al=2\al_1+4\al_3+6\al_4+3\al_2+5\al_5+4\al_6.\eqno(4.2)$$
Then
$$(\widehat\al,\al_r)=0\qquad\for\;\;r\in\ol{1,5}\eqno(4.3)$$
by the Dynkin diagram of $E_6$. Thanks to (2.9),
$$[\widehat\al,{\cal G}^{D_5}]=0.\eqno(4.4)$$
By Schur's Lemma, $\widehat\al|_V=c\sum_{i=1}^{16}x_i\ptl_{x_i}$.
According to the coefficients of $x_1\ptl_{x_1}$ in (3.1) with the
data in Table 2 and (4.2), we have
$$\widehat\al|_{\cal A}=\sum_{i=1}^{16}x_i\ptl_{x_i}=D,\eqno(4.5)$$
 the degree operator on ${\cal A}$.

Recall that the Lie bracket in the algebra $\mbb{A}$ (cf. (2.40)) is
given by the commutator
$$[d_1,d_2]=d_1d_2-d_2d_1.\eqno(4.6)$$
Set
$${\cal D}=\sum_{i=1}^{16}\mbb{C}\ptl_{x_i}.\eqno(4.7)$$
Then ${\cal D}$ forms an $o(10,\mbb{C})$-module with respect to the
action
$$B(d)=[B|_{\cal A},\ptl]\qquad\for\;\;B\in
o(10,\mbb{C}),\;\ptl\in{\cal D}.\eqno(4.8)$$ On the other hand,
${\cal G}_\pm$ (cf. (2.22)-(2.32)) form $o(10,\mbb{C})$-modules with
respect to the action
$$B(u)=[\nu(B),u]\qquad \for\;\;B\in
o(10,\mbb{C}),\;u\in{\cal G}_\pm.\eqno(4.9)$$ According to (2.36)
and (2.37), the linear map determined by $\eta_i\mapsto x_i$ for
$i\in\ol{1,16}$ gives an $o(10,\mbb{C})$-module isomorphism from
${\cal G}_-$ to $V$. Moreover, (2.12) and (2.34) implies that the
linear map determined by $\xi_i\mapsto \ptl_{x_i}$ for
$i\in\ol{1,16}$ gives an $o(10,\mbb{C})$-module isomorphism from
${\cal G}_+$ to ${\cal D}$. Hence we define the action of ${\cal
G}_+$ on ${\cal A}$ by
$$\xi_i|_{\cal
A}=\ptl_{x_i}\qquad\for\;\;i\in\ol{1,16}.\eqno(4.10)$$

Recall the Witt Lie subalgebra of $\mbb{A}$:
$${\cal W}_{16}=\sum_{i=1}^{16}{\cal A}\ptl_{x_i}.\eqno(4.11)$$
Now we want to find the differential operators
$P_1,P_2,...,P_{16}\in {\cal W}_{16}$ such that the following action
matches the structure of ${\cal G}^{E_6}$:
$$\eta_i|_{\cal A}=P_i\qquad\for\;\;i\in\ol{1,16}.\eqno(4.12)$$
Imposing
$$[\ptl_{x_1},P_1]=[E_{\al_6},E_{-\al_6}]|_{\cal A}=-\al_6|_{\cal
A}=2x_1\ptl_{x_1}+\sum_{i=2}^{10}x_i\ptl_{x_i}+x_{12}\ptl_{x_{12}},\eqno(4.13)$$
we take
$$P_1=x_1(\sum_{i=1}^{10}x_i\ptl_{x_i}+x_{12}\ptl_{x_{12}})+P_1',\eqno(4.14)$$
where $P_1'$ is a differential operator such that
$[\ptl_{x_1},P_1']=0$. Moreover,
$$[\ptl_{x_r},x_1(\sum_{i=1}^{10}x_i\ptl_{x_i}+x_{12}\ptl_{x_{12}})]=x_1\ptl_{x_r}\qquad\for\;\;r\in\{\ol{2,10},12\}.
\eqno(4.15)$$ Wanting $[\ptl_{x_r},P_1]\in {\cal G}^{D_5}|_{\cal
A}=o(10,\mbb{C})|_{\cal A}$ (cf. (2.49)-(2.71)), we take
\begin{eqnarray*} P_1&=&
x_1(\sum_{i=1}^{10}x_i\ptl_{x_i}+x_{12}\ptl_{x_{12}})-(x_2x_9-x_3x_6+x_4x_5)\ptl_{x_{11}}
\\ &&-(x_2x_{10}-x_3x_8+x_5x_7)\ptl_{x_{13}}-(x_2x_{12}-x_4x_8+x_6x_7)\ptl_{x_{14}}
\\ &&+(x_3x_{12}-x_4x_{10}+x_7x_9)\ptl_{x_{15}}-(x_5x_{12}-x_6x_{10}+x_8x_9)\ptl_{x_{16}}
\\&=&
x_1D-\zeta_1\ptl_{x_{11}}-\zeta_2\ptl_{x_{13}}-\zeta_3\ptl_{x_{14}}-\zeta_4\ptl_{x_{15}}-\zeta_{10}\ptl_{x_{16}}
\hspace{4.8cm}(4.16)\end{eqnarray*} by (2.52), (2.56), (2.59),
(2.61), (2.62), (2.64)-(2.68), (3.8), (3.14)-(3.16), (3.18) and
(4.5). Then
$$[\ptl_{x_s},P_1]=[\xi_s,\eta_1]|_{\cal
A}\qquad\for\;\;s\in\ol{1,16}\eqno(4.17)$$ due to (2.42)-(2.48).

Since $[E_{-\al_5},\eta_1]=\eta_2$ by (2.10), we take
\begin{eqnarray*}\qquad P_2&=&[E_{-\al_5}|_{\cal A},\eta_1|_{\cal A}]=-[(E_{5,4}-E_{9,10})|_{\cal
A},P_1]\\&=&
x_2D-\zeta_1\ptl_{x_9}-\zeta_2\ptl_{x_{10}}-\zeta_3\ptl_{x_{12}}+\zeta_5\ptl_{x_{15}}-\zeta_9\ptl_{x_{16}}
\hspace{4.2cm}(4.18)\end{eqnarray*} by (2.44) and (3.12). Note that
 $[E_{-\al_4},\eta_2]=\eta_3$ by (2.10).
 Hence (3.11) gives
\begin{eqnarray*}\qquad P_3&=&[E_{-\al_4}|_{\cal A},\eta_2|_{\cal A}]=-[(E_{4,3}-E_{8,9})|_{\cal
A},P_2]\\&=&
x_3D+\zeta_1\ptl_{x_6}+\zeta_2\ptl_{x_8}+\zeta_4\ptl_{x_{12}}+\zeta_5\ptl_{x_{14}}-\zeta_8\ptl_{x_{16}}.
\hspace{4.3cm}(4.19)\end{eqnarray*} Thanks to
$-[E_{-\al_3},\eta_3]=\eta_4$, (2.50) and (3.10), we have
\begin{eqnarray*} \qquad P_4&=&-[E_{-\al_3}|_{\cal A},\eta_3|_{\cal A}]=[(E_{3,2}-E_{7,9})|_{\cal
A},P_3]\\
&=&
x_4D-\zeta_1\ptl_{x_5}+\zeta_3\ptl_{x_8}-\zeta_4\ptl_{x_{10}}-\zeta_5\ptl_{x_{13}}+\zeta_7\ptl_{x_{16}}.
\hspace{4.3cm}(4.20)\end{eqnarray*} Observe that
$[E_{-\al_2},\eta_3]=\eta_5$ by (2.10).  So (3.13) yields
\begin{eqnarray*} \qquad P_5&=&[E_{-\al_2}|_{\cal A},\eta_3|_{\cal A}]=-[(E_{10,4}-E_{9,5})|_{\cal
A},P_2]\\
&=&
x_5D-\zeta_1\ptl_{x_4}-\zeta_2\ptl_{x_7}-\zeta_{10}\ptl_{x_{12}}+\zeta_9\ptl_{x_{14}}-\zeta_8\ptl_{x_{15}}.
\hspace{4.2cm}(4.21)\end{eqnarray*}

Since $[E_{-\al_2},\eta_4]=\eta_6$, (2.13) implies
\begin{eqnarray*} \qquad P_6&=&[E_{-\al_2}|_{\cal A},\eta_4|_{\cal A}]=-[(E_{10,4}-E_{9,5})|_{\cal
A},P_4]\\
&=&
x_6D+\zeta_1\ptl_{x_3}-\zeta_3\ptl_{x_7}+\zeta_{10}\ptl_{x_{10}}-\zeta_9\ptl_{x_{13}}+\zeta_7\ptl_{x_{15}}.
\hspace{4.2cm}(4.22)\end{eqnarray*} As
$-[E_{-\al_1},\eta_4]=\eta_7$,
 we get by (3.9) that
\begin{eqnarray*} \qquad P_7&=&-[E_{-\al_1}|_{\cal A},\eta_4|_{\cal A}]=[(E_{2,1}-E_{6,7})|_{\cal
A},P_4]\\ &=&
x_7D-\zeta_2\ptl_{x_5}-\zeta_3\ptl_{x_6}+\zeta_4\ptl_{x_9}+\zeta_5\ptl_{x_{11}}-\zeta_6\ptl_{x_{16}}.
\hspace{4.5cm}(4.23)\end{eqnarray*} Thanks to
$[E_{-\al_2},\eta_7]=\eta_8$, (3.13) gives
\begin{eqnarray*}\qquad P_8&=&[E_{-\al_2}|_{\cal A},\eta_7|_{\cal A}=-[(E_{10,4}-E_{9,5})|_{\cal
A},P_7]
\\ &=&
x_8D+\zeta_2\ptl_{x_3}+\zeta_3\ptl_{x_4}-\zeta_{10}\ptl_{x_9}+\zeta_9\ptl_{x_{11}}-\zeta_6\ptl_{x_{15}}.
\hspace{4.4cm}(4.24)\end{eqnarray*}

The fact $[E_{-\al_4},\eta_6]=\eta_9$ yields
\begin{eqnarray*}\qquad P_9&=&[E_{-\al_4}|_{\cal A},\eta_6|_{\cal A}]=-[(E_{4,3}-E_{8,9})|_{\cal
A},P_6]
\\ &=&
x_9D-\zeta_1\ptl_{x_2}+\zeta_4\ptl_{x_7}-\zeta_{10}\ptl_{x_8}-\zeta_8\ptl_{x_{13}}+\zeta_7\ptl_{x_{14}}.
\hspace{4.4cm}(4.25)\end{eqnarray*}
 by (3.11). As
$[E_{-\al_4},\eta_8]=\eta_{10}$, we find
\begin{eqnarray*} \qquad P_{10}&=&[E_{-\al_4}|_{\cal A},\eta_8|_{\cal A}=-[(E_{4,3}-E_{8,9})|_{\cal
A},P_8]
\\ &=&
x_{10}D-\zeta_2\ptl_{x_2}-\zeta_4\ptl_{x_4}+\zeta_{10}\ptl_{x_6}+\zeta_8\ptl_{x_{11}}-\zeta_6\ptl_{x_{14}}.
\hspace{4.1cm}(4.26)\end{eqnarray*} by (3.11). Moreover, the fact
$-[E_{-\al_5},\eta_9]=\eta_{11}$ implies
\begin{eqnarray*} \qquad P_{11}&=&-[E_{-\al_5}|_{\cal A},\eta_9|_{\cal A}]=[(E_{5,4}-E_{9,10})|_{\cal
A},P_9]
\\ &=&
x_{11}D-\zeta_1\ptl_{x_1}+\zeta_5\ptl_{x_7}+\zeta_9\ptl_{x_8}+\zeta_8\ptl_{x_{10}}-\zeta_7\ptl_{x_{12}}.
\hspace{4.3cm}(4.27)\end{eqnarray*}
 by (3.12).
  Since $-[E_{-\al_3},\eta_{10}]=\eta_{12}$,
\begin{eqnarray*}\qquad P_{12}&=&-[E_{-\al_3}|_{\cal A},\eta_{10}|_{\cal A}=[(E_{3,2}-E_{7,8})|_{\cal
A},P_{10}]
\\ &=&
x_{12}D-\zeta_3\ptl_{x_2}+\zeta_4\ptl_{x_3}-\zeta_{10}\ptl_{x_5}-\zeta_7\ptl_{x_{11}}+\zeta_6\ptl_{x_{13}}
\hspace{4.1cm}(4.28)\end{eqnarray*} by (3.10).

Observing $-[E_{-\al_1},\eta_{11}]=\eta_{13}$, we have
\begin{eqnarray*}\qquad P_{13}&=&-[E_{-\al_1}|_{\cal A},\eta_{11}|_{\cal A}]=[(E_{2,1}-E_{6,7})|_{\cal
A},P_{11}]
\\ &=&
x_{13}D-\zeta_2\ptl_{x_1}-\zeta_5\ptl_{x_4}-\zeta_9\ptl_{x_6}-\zeta_8\ptl_{x_9}+\zeta_6\ptl_{x_{12}}.
\hspace{4.3cm}(4.29)\end{eqnarray*} by (3.9). The fact
$-[E_{-\al_3},\eta_{13}]=\eta_{14}$ gives
\begin{eqnarray*}\qquad P_{14}&=&-[E_{-\al_3}|_{\cal A},\eta_{13}|_{\cal A}]=[(E_{3,2}-E_{7,8})|_{\cal
A},P_{13}]
\\ &=&
x_{14}D-\zeta_3\ptl_{x_1}+\zeta_5\ptl_{x_3}+\zeta_9\ptl_{x_5}+\zeta_7\ptl_{x_9}-\zeta_6\ptl_{x_{10}}.
\hspace{4.3cm}(4.30)\end{eqnarray*}
 by (3.10). As
$-[E_{-\al_4},\eta_{14}]=\eta_{15}$, we get
\begin{eqnarray*}\qquad P_{15}&=&-[E_{-\al_4}|_{\cal A},\eta_{14}|_{\cal A}]=[(E_{4,3}-E_{8,9})|_{\cal
A},P_{14}]
\\ &=&
x_{15}D-\zeta_4\ptl_{x_1}+\zeta_5\ptl_{x_2}-\zeta_8\ptl_{x_5}+\zeta_7\ptl_{x_6}-\zeta_6\ptl_{x_8}.
\hspace{4.4cm}(4.31)\end{eqnarray*}
 by (3.11).
 Since
$-[E_{-\al_2},\eta_{15}]=\eta_{16}$, we have
\begin{eqnarray*} \qquad P_{16}&=&-[E_{-\al_2}|_{\cal A},\eta_{15}|_{\cal A}]=[(E_{10,4}-E_{9,5})|_{\cal
A},P_{15}]\\ &=&
x_{16}D-\zeta_{10}\ptl_{x_1}-\zeta_9\ptl_{x_2}-\zeta_8\ptl_{x_3}+\zeta_7\ptl_{x_4}-\zeta_6\ptl_{x_7}.
\hspace{4.3cm}(4.32)\end{eqnarray*}
 by (3.13).

Set
$${\cal P}=\sum_{i=1}^{16}\mbb{C}P_i,\qquad{\cal C}_0=o(10,\mbb{C})|_{\cal
A}+\mbb{C}D\eqno(4.33)$$ (cf. (2.49)-(2.71) and (4.5)) and
$${\cal
C}={\cal P}+{\cal C}_0+{\cal D}\eqno(4.34)$$ (cf. (4.7)). The we
have:\psp

{\bf Theorem 4.1}. {\it The space ${\cal C}$ forms a Lie subalgebra
of the Witt algebra ${\cal W}_{16}$ (cf. (4.11)). Moreover, the
linear map $\vt$ determined by
$$\vt(\xi_i)=\ptl_{x_i},\;\;\vt(\eta_i)=P_i,\;\;\vt(u)=\nu^{-1}(u)|_{\cal
A}\qquad\for\;\;i\in\ol{1,16},\;u\in{\cal G}^{D_5}\eqno(4.35)$$ (cf.
(2.42)-(2.47)) and
$$\vt(\al_6)=-2x_1\ptl_{x_1}-\sum_{i=2}^{10}x_i\ptl_{x_i}-x_{12}\ptl_{x_{12}}\eqno(4.36)$$
(cf. (4.1)) gives a Lie algebra isomorphism from ${\cal G}^{E_6}$ to
${\cal C}$.}

{\it Proof}. Since ${\cal D}\cong {\cal G}_+$ as ${\cal
G}^{D_5}$-modules, we have
$${\cal G}_0+{\cal G}_+\stl{\vt}{\cong}{\cal C}_0+{\cal D}\eqno(4.37)$$ as Lie algebras. Denote by
$U({\cal G})$ the universal enveloping algebra of a Lie algebra
${\cal G}$. Note that
$${\cal B}_-={\cal G}_0+{\cal G}_-,\qquad{\cal B}_+={\cal G}_0+{\cal G}_+\eqno(4.38)$$
are also  Lie subalgebras of ${\cal G}^{E_6}$ and
$${\cal G}^{E_6}={\cal B}_-\oplus {\cal G}_+={\cal G}_-\oplus {\cal B}_+.\eqno(4.39)$$
We define a one-dimensional ${\cal B}_-$-module $\mbb{C}u_0$ by
$$w(u_0)=0\qquad\for\;\;w\in{\cal B}_-\bigcap{\cal G}^{D_5},\;\;\widehat\al(u_0)=-16u_0\eqno(4.40)$$
(cf. (4.2)). Let
$$\Psi=U({\cal G}^{E_6})\otimes_{{\cal B}_-}\mbb{C}u_0\cong U({\cal G}_+)\otimes_\mbb{C} \mbb{C}u_0\eqno(4.41)$$
be the induced ${\cal G}^{E_6}$-module.

Reacll that $\mbb{N}$ is the set of nonnegative integers. Let
$${\cal
A}'=\mbb{C}[\ptl_{x_1},\ptl_{x_2},...,\ptl_{x_{16}}].\eqno(4.42)$$
We define an action of the associative algebra $\mbb{A}$ (cf.
(2.40)) on ${\cal A}'$ by
$$\ptl_{x_i}(\prod_{j=1}^{16}\ptl_{x_j}^{\be_j})=\ptl_{x_i}^{\be_i+1}\prod_{i\neq
j\in\ol{1,16}}\prod_{j=1}^{16}\ptl_{x_j}^{\be_j}\eqno(4.43)$$ and
$$x_i(\prod_{j=1}^{16}\ptl_{x_j}^{\be_j})=-\be_i\ptl_{x_i}^{\be_i-1}\prod_{i\neq
j\in\ol{1,16}}\prod_{j=1}^{16}\ptl_{x_j}^{\be_j}\eqno(4.44)$$ for
$i\in\ol{1,16}$. Since
$$[-x_i,\ptl_{x_j}]=[\ptl_{x_i},x_j]=\dlt_{i,j}\qquad\for\;\;i,j\in\ol{1,16},\eqno(4.45)$$
the above action gives an associative algebra representation of
$\mbb{A}$. Thus it also gives a Lie algebra representation of
$\mbb{A}$ (cf. (4.6)).
 It is straightforward to verify that
$$[d|_{{\cal A}'},\ptl|_{{\cal A}'}]=[d,\ptl]|_{{\cal
A}'}\qquad\for\;\;d\in{\cal C}_0,\;\ptl\in{\cal D}.\eqno(4.46)$$

Define linear map $\vs: \Psi\rta {\cal A}'$ by
$$\vs(\prod_{i=1}^{16}\xi_i^{\ell_i}\otimes
u_0)=\prod_{i=1}^{16}\ptl_{x_i}^{\ell_i}\qquad(\ell_1,...,\ell_{16})\in\mbb{N}^{16}.\eqno(4.47)$$
According to (2.49)-(2.71), (4.43) and (4.44),
$$D(1)=-16,\;\;d(1)=0\qquad\for\;\;d\in o(10,\mbb{C})|_{\cal A}.\eqno(4.48)$$
Moreover, (4.40), (4.41),  (4.43), (4.44) and (4.48) imply
$$\vs(\xi(v))=\vt(\xi)\vs(v)\qquad\for\;\;\xi\in{\cal
G}_0,\;v\in\Psi.\eqno(4.49)$$

Now (4.41) and (4.43) imply
$$\vs(w(u))=\vt(w)(\vs(u))\qquad\for\;\;w\in{\cal
B}_+,\;\;u\in\Psi.\eqno(4.50)$$ Thus  we have
$$\vs w|_{\Psi}\vs^{-1}=\vt(w)|_{{\cal A}'}\qquad\for\;\;w\in{\cal
B}_+.\eqno(4.51)$$ On the other hand, the linear map
$$\psi(v)=\vs v|_{\Psi}\vs^{-1}\qquad\for\;\; v\in{\cal
G}^{E_6}\eqno(4.52)$$ is a Lie algebra monomorphism from ${\cal
G}^{E_6}$ to $\mbb{A}|_{{\cal A}'}$. According to (4.17) and (4.45),
$$\psi(\eta_1)=P_1|_{{\cal A}'}.\eqno(4.53)$$
By the constructions of $P_2,...,P_{16}$ in (4.18)-(4.32), we have
$$\psi(\eta_i)=P_i|_{{\cal
A}'}\qquad\for\;\;i\in\ol{2,16}.\eqno(4.54)$$ Therefore, we have
$$\psi(v)=\vt(v)|_{{\cal A}'}\qquad\for\;\;v\in{\cal
G}^{E_6}.\eqno(4.55)$$ In particular, ${\cal C}|_{{\cal
A}'}=\vt({\cal G}^{E_6})|_{{\cal A}'}=\psi({\cal G}^{E_6})$ forms a
Lie algebra. Since the linear map $d\mapsto d|_{{\cal A}'}$ for
$d\in {\cal C}$ is injective, we have that ${\cal C}$ forms a Lie
subalgebra of $\mbb{A}$ and $\vt$ is a Lie algebra
isomorphism.$\qquad\Box$ \psp

By the above theorem, a Lie group of type $E_6$ is generated by the
real spinor transformations corresponding to (2.49)-(2.71), the real
translations and  dilations in $\sum_{i=1}^{16}\mbb{R}x_i$, and the
following fractional transformations:
$$\wp_{1b}(x_i)=\frac{x_i}{1-bx_1},\;\;i\in\{\ol{1,10},12\},\;\;\wp_{1b}(x_{11})
=x_{11}-\frac{b(x_2x_9-x_3x_6+x_4x_5)}{1-bx_1},\eqno(4.56)$$
$$ \wp_{1b}(x_{13})=x_{13}-\frac{b(x_2x_{10}-x_3x_8+x_5x_7)}{1-bx_1},\;\wp_{1b}(x_{14})=x_{14}-\frac{
b(x_2x_{12}-x_4x_8+x_6x_7)}{1-bx_1},\eqno(4.57)$$
$$
\wp_{1b}(x_{15})=x_{15}+\frac{b(x_3x_{12}-x_4x_{10}+x_7x_9)}{1-bx_1},\eqno(4.58)$$
$$\wp_{1b}(x_{16})=x_{16}-\frac{b(x_5x_{12}-x_6x_{10}+x_8x_9)}{1-bx_1};\eqno(4.59)$$
$$\wp_{2b}(x_i)=\frac{x_i}{1-bx_2},\;\;i\in\ol{1,14}\setminus\{9,10,12\},\;\;\wp_{2
b}(x_9)=x_9-\frac{b(x_1x_{11}-x_3x_6+x_4x_5)}{1-bx_2},\eqno(4.60)$$
$$ \wp_{2b}(x_{10})=x_{10}-\frac{b(x_1x_{13}-x_3x_8+x_5x_7)}{1-bx_2},\;
\wp_{2b}(x_{12})=x_{12}-\frac{b(x_1x_{14}-x_4x_8+x_6x_7)}{1-bx_2},\eqno(4.61)$$
$$
\wp_{2b}(x_{15})=x_{15}-\frac{b(x_3x_{14}-x_4x_{13}+x_7x_{11})}{1-bx_2},\eqno(4.62)$$
$$\wp_{2b}(x_{16})=x_{16}+\frac{b(x_5x_{14}-x_6x_{13}+x_8x_{11})}{1-bx_2};\eqno(4.63)$$
$$\wp_{3b}(x_i)=\frac{x_i}{1-bx_3},\;i\in\ol{1,15}\setminus\{6,8,12,14\},\;\wp_{3
b}(x_6)=x_6+\frac{b(x_1x_{11}+x_2x_9+x_4x_5)}{1-bx_3},\eqno(4.64)$$
$$ \wp_{3b}(x_8)=x_8+\frac{b(x_1x_{13}+x_2x_{10}+x_5x_7)}{1-bx_3}
,\;
\wp_{3b}(x_{12})=x_{12}+\frac{b(x_1x_{15}+x_4x_{10}-x_7x_9)}{1-bx_3},\eqno(4.65)$$
$$
\wp_{3b}(x_{14})=x_{14}-\frac{b(x_2x_{15}-x_4x_{13}+x_7x_{11})}{1-bx_3},\eqno(4.66)$$
$$\wp_{3b}(x_{16})=x_{16}-\frac{b(x_5x_{15}+x_9x_{13}-x_{10}x_{11})}{1-bx_3};\eqno(4.67)$$
$$\wp_{4b}(x_i)=\frac{x_i}{1-bx_4},\;i\in\ol{1,15}\setminus\{5,8,10,13\},\;\wp_{4
b}(x_5)= x_5-\frac{b(x_1x_{11}+x_2x_9-x_3x_6)}{1-bx_4},\eqno(4.68)$$
$$ \wp_{4b}(x_8)=x_8+\frac{b(x_1x_{14}+x_2x_{12}+x_6x_7)}{1-bx_4}
,\;
\wp_{4b}(x_{10})=x_{10}-\frac{b(x_1x_{15}-x_3x_{12}-x_7x_9)}{1-bx_4},\eqno(4.69)$$
$$
\wp_{4b}(x_{13})=x_{13}+\frac{b(x_2x_{15}+x_3x_{14}+x_7x_{11})}{1-bx_4},\eqno(4.70)$$
$$\wp_{4b}(x_{16})=x_{16}-\frac{b(x_6x_{15}+x_9x_{14}-x_{11}x_{12})}{1-bx_4};\eqno(4.71)$$
$$\wp_{5b}(x_i)=\frac{x_i}{1-bx_5},\;\;\wp_{5
b}(x_4)= x_4-\frac{b(x_1x_{11}+x_2x_9-x_3x_6)}{1-bx_5},\eqno(4.72)$$
$i\in\ol{1,16}\setminus\{4,7,12,14,15\}$,
$$ \wp_{5b}(x_7)=x_7-\frac{b(x_1x_{13}+x_2x_{10}-x_3x_8)}{1-bx_5}
,\;
\wp_{5b}(x_{12})=x_{12}-\frac{b(x_1x_{16}-x_6x_{10}+x_8x_9)}{1-bx_5},\eqno(4.73)$$
$$ \wp_{5b}(x_{14})=x_{14}+\frac{b(x_2x_{16}+x_6x_{13}-x_8x_{11})}{1-bx_5}
,\eqno(4.74)$$
$$\wp_{5b}(x_{15})=x_{15}-\frac{b(x_3x_{16}+x_9x_{13}-x_{10}x_{11})}{1-bx_5};\eqno(4.75)$$
$$\wp_{6b}(x_i)=\frac{x_i}{1-bx_6},\;\;\wp_{6
b}(x_3)= x_3+\frac{b(x_1x_{11}+x_2x_9+x_4x_5)}{1-bx_6},\eqno(4.76)$$
$i\in\ol{1,16}\setminus\{3,7,10,13,15\}$,
$$ \wp_{6b}(x_7)=x_7-\frac{b(x_1x_{14}+x_2x_{12}-x_4x_8)}{1-bx_6}
,\;
\wp_{6b}(x_{10})=x_{10}+\frac{b(x_1x_{16}+x_5x_{12}+x_8x_9)}{1-bx_6},\eqno(4.77)$$
$$ \wp_{6b}(x_{13})=x_{13}-\frac{b(x_2x_{16}-x_5x_{14}-x_8x_{11})}{1-bx_6}
,\eqno(4.78)$$
$$\wp_{6b}(x_{15})=x_{15}-\frac{b(x_4x_{16}+x_9x_{14}-x_{11}x_{12})}{1-bx_6};\eqno(4.79)$$
$$\wp_{7b}(x_i)=\frac{x_i}{1-bx_7},\;i\in\ol{1,15}\setminus\{5,6,9,11\},\;\wp_{7
b}(x_5)=
x_5-\frac{b(x_1x_{13}+x_2x_{10}-x_3x_8)}{1-bx_7},\eqno(4.80)$$
$$ \wp_{7b}(x_6)=x_6-\frac{b(x_1x_{14}+x_2x_{12}-x_4x_8)}{1-bx_7}
,\;
\wp_{7b}(x_9)=x_9+\frac{b(x_1x_{15}-x_3x_{12}+x_4x_{10})}{1-bx_7},\eqno(4.81)$$
$$ \wp_{7b}(x_{11})=x_{11}-\frac{b(x_2x_{15}+x_3x_{14}-x_4x_{13})}{1-bx_7}
,\eqno(4.82)$$
$$\wp_{7b}(x_{16})=x_{16}-\frac{b(x_8x_{15}+x_{10}x_{14}
 -x_{12}x_{13})}{1-bx_7};\eqno(4.83)$$
$$\wp_{8b}(x_i)=\frac{x_i}{1-bx_8},\;\;\wp_{8
b}(x_3)= x_3+\frac{b(x_1x_{13}-x_3x_8+x_5x_7)}{1-bx_8},\eqno(4.84)$$
$i\in\ol{1,16}\setminus\{3,4,9,11,15\}$,
$$ \wp_{8b}(x_4)=x_4+\frac{b(x_1x_{14}+x_2x_{12}+x_6x_7)}{1-bx_8}
,\;
\wp_{8b}(x_9)=x_9-\frac{b(x_1x_{16}+x_5x_{12}-x_6x_{10})}{1-bx_8},\eqno(4.85)$$
$$ \wp_{8b}(x_{11})=x_{11}+\frac{b(x_2x_{16}-x_5x_{14}+x_6x_{13})}{1-bx_8}
,\eqno(4.86)$$
$$\wp_{8b}(x_{15})=x_{15}-\frac{b(x_7x_{16}+x_{10}x_{14}
 -x_{12}x_{13})}{1-bx_8};\eqno(4.87)$$
$$\wp_{9b}(x_i)=\frac{x_i}{1-bx_9},\;\;\wp_{9
b}(x_2)= x_2-\frac{b(x_1x_{11}-x_3x_6+x_4x_5)}{1-bx_9},\eqno(4.88)$$
$i\in\ol{1,16}\setminus\{2,7,8,13,14\}$,
$$ \wp_{9b}(x_7)=x_7+\frac{b(x_1x_{15}-x_3x_{12}+x_4x_{10})}{1-bx_9}
,\;
\wp_{9b}(x_8)=x_8-\frac{b(x_1x_{16}+x_5x_{12}-x_6x_{10})}{1-bx_9},\eqno(4.89)$$
$$ \wp_{9b}(x_{13})=x_{13}-\frac{b(x_3x_{16}+x_5x_{15}-x_{10}x_{11})}{1-bx_9}
,\eqno(4.90)$$
$$\wp_{9b}(x_{14})=x_{14}-\frac{b(x_4x_{16}+x_6x_{15}
 -x_{11}x_{12})}{1-bx_9};\eqno(4.91)$$
$$\wp_{10b}(x_i)=\frac{x_i}{1-bx_{10}},\;\;\wp_{10
b}(x_2)=
x_2-\frac{b(x_1x_{13}-x_3x_8+x_5x_7)}{1-bx_{10}},\eqno(4.92)$$
$i\in\ol{1,16}\setminus\{2,4,6,11,14\}$,
$$ \wp_{10b}(x_4)=x_4-\frac{b(x_1x_{15}-x_3x_{12}-x_7x_9)}{1-bx_{10}}
,\eqno(4.93)$$
$$\wp_{10b}(x_6)=x_6+\frac{b(x_1x_{16}+x_5x_{12}+x_8x_9)}{1-bx_{10}},\eqno(4.94)$$
$$ \wp_{10b}(x_{11})=x_{11}+\frac{b(x_3x_{16}+x_5x_{15}+x_9x_{13})}{1-bx_{11}}
,\eqno(4.95)$$
$$\wp_{10b}(x_{14})=x_{14}-\frac{b(x_7x_{16}+x_8x_{15}
 -x_{12}x_{13})}{1-bx_{10}};\eqno(4.96)$$
$$\wp_{11b}(x_i)=\frac{x_i}{1-bx_{11}},\;\;\wp_{11
b}(x_1)= x_1-\frac{b(x_2x_9-x_3x_6+x_4x_5)}{1-bx_{11}},\eqno(4.97)$$
$i\in\ol{1,16}\setminus\{1,7,8,10,12\}$,
$$ \wp_{11b}(x_7)=x_7-\frac{b(x_2x_{15}+x_3x_{14}-x_4x_{13})}{1-bx_{11}}
,\eqno(4.98)$$
$$\wp_{11b}(x_8)=x_8+\frac{b(x_2x_{16}-x_5x_{14}+x_6x_{13})}{1-bx_{11}},\eqno(4.99)$$
$$ \wp_{11b}(x_{10})=x_{10}+\frac{b(x_3x_{16}+x_5x_{15}+x_9x_{13})}{1-bx_{11}}
,\eqno(4.100)$$
$$\wp_{11b}(x_{12})=x_{12}+\frac{b(x_4x_{16}+x_6x_{15}
 +x_9x_{14})}{1-bx_{11}};\eqno(4.101)$$
$$\wp_{12b}(x_i)=\frac{x_i}{1-bx_{12}},\;\;\wp_{12
b}(x_2)=
x_2-\frac{b(x_1x_{14}+x_6x_7-x_4x_8)}{1-bx_{12}},\eqno(4.102)$$
$i\in\ol{1,16}\setminus\{2,3,5,11,13\}$,
$$ \wp_{12b}(x_3)=x_3+\frac{b(x_1x_{15}-x_7x_9+x_4x_{10})}{1-bx_{12}}
,\eqno(4.103)$$
$$\wp_{12b}(x_5)=x_5-\frac{b(x_1x_{16}+x_8x_9-x_6x_{10})}{1-bx_{12}},\eqno(4.104)$$
$$ \wp_{12b}(x_{11})=x_{11}+\frac{b(x_4x_{16}+x_6x_{15}+x_9x_{14})}{1-bx_{12}}
,\eqno(4.105)$$
$$\wp_{12b}(x_{13})=x_{13}+\frac{b(x_7x_{16}+x_8x_{15}
 +x_{10}x_{14})}{1-bx_{12}};\eqno(4.106)$$
$$\wp_{13b}(x_i)=\frac{x_i}{1-bx_{13}},\;\;\wp_{13
b}(x_1)=
x_1-\frac{b(x_1x_{13}-x_3x_8+x_5x_7)}{1-bx_{13}},\eqno(4.107)$$
$i\in\ol{1,16}\setminus\{1,4,6,9,12\}$,
$$ \wp_{13b}(x_4)=x_4+\frac{b(x_2x_{15}+x_3x_{14}+x_7x_{11})}{1-bx_{13}}
,\eqno(4.108)$$
$$\wp_{13b}(x_6)=x_6-\frac{b(x_2x_{16}-x_5x_{14}-x_8x_{11})}{1-bx_{13}},\eqno(4.109)$$
$$ \wp_{13b}(x_9)=x_9-\frac{b(x_3x_{16}+x_5x_{15}-x_{10}x_{11})}{1-bx_{13}}
,\eqno(4.110)$$
$$\wp_{13b}(x_{12})=x_{12}+\frac{b(x_7x_{16}+x_8x_{15}
 +x_{10}x_{14})}{1-bx_{13}};\eqno(4.111)$$
$$\wp_{14b}(x_i)=\frac{x_i}{1-bx_{14}},\;\;\wp_{14
b}(x_1)=
x_1-\frac{b(x_2x_{12}-x_4x_8+x_6x_7)}{1-bx_{14}},\eqno(4.112)$$
$i\in\ol{1,16}\setminus\{1,3,5,9,10\}$,
$$ \wp_{14b}(x_3)=x_3-\frac{b(x_2x_{15}-x_4x_{13}+x_7x_{11})}{1-bx_{14}}
,\eqno(4.113)$$
$$\wp_{14b}(x_5)=x_5+\frac{b(x_2x_{16}+x_6x_{13}-x_8x_{11})}{1-bx_{14}},\eqno(4.114)$$
$$ \wp_{14b}(x_9)=x_9-\frac{b(x_4x_{16}+x_6x_{15}-x_{11}x_{12})}{1-bx_{14}}
,\eqno(4.115)$$
$$\wp_{14b}(x_{10})=x_{10}-\frac{b(x_7x_{16}+x_8x_{15}
 -x_{12}x_{13})}{1-bx_{14}};\eqno(4.116)$$
$$\wp_{15b}(x_i)=\frac{x_i}{1-bx_{15}},\;\;\wp_{15
b}(x_1)=
x_1+\frac{b(x_3x_{12}-x_4x_{10}+x_7x_9)}{1-bx_{15}},\eqno(4.117)$$
$i\in\ol{1,16}\setminus\{1,2,5,6,8\}$,
$$ \wp_{15b}(x_2)=x_2-\frac{b(x_3x_{14}-x_4x_{13}+x_7x_{11})}{1-bx_{15}}
,\eqno(4.118)$$
$$\wp_{15b}(x_5)=x_5-\frac{b(x_3x_{16}+x_9x_{13}-x_{10}x_{11})}{1-bx_{15}},\eqno(4.119)$$
$$ \wp_{15b}(x_6)=x_6-\frac{b(x_4x_{16}+x_9x_{14}-x_{11}x_{12})}{1-bx_{15}}
,\eqno(4.120)$$
$$\wp_{15b}(x_8)=x_8-\frac{b(x_7x_{16}+x_{10}x_{14}
 -x_{12}x_{13})}{1-bx_{15}};\eqno(4.121)$$

$$\wp_{16b}(x_i)=\frac{x_i}{1-bx_{16}},\;\;\wp_{16
b}(x_1)=
x_1-\frac{b(x_5x_{12}-x_6x_{10}+x_8x_9)}{1-bx_{16}},\eqno(4.122)$$
$i\in\ol{1,16}\setminus\{1,2,3,4,7\}$,
$$ \wp_{16b}(x_2)=x_2+\frac{b(x_5x_{14}-x_6x_{13}+x_8x_{11})}{1-bx_{16}}
,\eqno(4.123)$$
$$\wp_{16b}(x_3)=x_3-\frac{b(x_5x_{15}+x_9x_{13}-x_{10}x_{11}))}{1-bx_{16}},\eqno(4.124)$$
$$ \wp_{16b}(x_4)=x_4-\frac{b(x_6x_{15}+x_9x_{14}-x_{11}x_{12})}{1-bx_{16}}
,\eqno(4.125)$$
$$\wp_{16b}(x_7)=x_7-\frac{b(x_8x_{15}+x_{10}x_{14}
 -x_{12}x_{13})}{1-bx_{16}};\eqno(4.126)$$
where $b\in\mbb{R}$.

\section{Functor from $D_5$-Mod to $E_6$-Mod}

 In this section, we construct a new functor from $D_5$-{\bf Mod} to $E_6$-{\bf
 Mod}.

Note that
\begin{eqnarray*}o(10,{\cal A})&=&\sum_{1\leq
p<q\leq 5}[{\cal A}(E_{p,n+q}-E_{q,n+p})+{\cal A}(E_{n+p,q}-E_{n+q,p})]\\
& &+\sum_{i,j=1}^5{\cal
A}(E_{i,j}-E_{n+j,n+i})\hspace{8cm}(5.1)\end{eqnarray*} forms a Lie
subalgebra of the matrix algebra $gl(10,{\cal A})$ over ${\cal A}$
with respect to the commutator, i.e.
$$[fB_1,gB_2]=fg[B_1,B_2]\qquad\for\;\;f,g\in{\cal A},\;B_1,B_2\in
gl(10,\mbb{C}).\eqno(5.2)$$ Moreover, we define the Lie algebra
$${\cal K}=o(10,{\cal A})\oplus {\cal A}\kappa\eqno(5.3)$$
with the Lie bracket:
$$[\xi_1+f\kappa,\xi_2+g\kappa]=[\xi_1,\xi_2]\qquad\for\;\;\xi_1,\xi_2\in
o(10,{\cal A}),\;f,g\in{\cal A}.\eqno(5.4)$$ Similarly, $gl(16,{\cal
A})$ becomes a Lie algebra with the Lie bracket as that in (5.2).
Recall the Witt algebra ${\cal W}_{16}=\sum_{i=1}^{16}{\cal
A}\ptl_{x_i}$, and Shen [Sg1-3] found a monomorphism $\Im$ from the
Lie algebra ${\cal W}_{16}$ to the Lie algebra of semi-product
${\cal W}_{16}+gl(16,{\cal A})$ defined by
$$\Im(\sum_{i=1}^{16}f_i\ptl_{x_i})=\sum_{i=1}^{16}f_i\ptl_{x_i}+\Im_1(\sum_{i=1}^{16}f_i\ptl_{x_i}),\;\;
\Im_1(\sum_{i=1}^{16}f_i\ptl_{x_i})=\sum_{i,j=1}^{16}\ptl_{x_i}(f_j)E_{i,j}.
\eqno(5.5)$$ According to our construction of $P_1$-$P_{16}$ in
(4.12)-(4.32),
$$\Im_1(P_i)=\sum_{r=1}^{16}x_r\Im_1([\xi_r,\eta_i]|_{\cal A})\qquad\for\;\;i\in\ol{1,16}.\eqno(5.6)$$

On the other hand,
$$\widehat{\cal K}={\cal W}_{16}\oplus {\cal K}\eqno(5.7)$$
becomes a Lie algebra with the Lie bracket
\begin{eqnarray*}&&[d_1+f_1B_1+f_2\kappa,d_2+g_1B_2+g_2\kappa]\\
&=&[d_1,d_2]+f_1g_1[B_1,B_2]+d_1(g_2)B_2
-d_2(g_1)B_1+(d_1(g_2)-d_2(g_1))\kappa\hspace{2.5cm}(5.8)\end{eqnarray*}
for $f_1,f_2,g_1,g_2\in{\cal A},\;B_1,B_2\in o(10,\mbb{C})$ and
$d_1,d_2\in{\cal W}_{16}$. Note
$${\cal G}_0={\cal
G}^{D_5}\oplus\mbb{C}\widehat\al\eqno(5.9)$$ by (2.32) and (4.2). So
there exists a Lie algebra monomorphism $\varrho:{\cal G}_0\rta
{\cal K}$ determined by
$$\varrho(\widehat\al)=2\kappa,\;\;\varrho(u)=\nu^{-1}(u)\qquad\for\;\;u\in{\cal
G}^{D_5}.\eqno(5.10)$$ Since $\Im$ is a Lie algebra monomorphism,
our construction of $P_1$-$P_{16}$ in (4.12)-(4.32) and (5.6) shows
that we have a Lie algebra monomorphism $\iota: {\cal G}^{E_6}\rta
\widehat{\cal K}$ given by
$$\iota(u)=u|_{\cal A}+\varrho(u)\qquad\for\;\;u\in{\cal
G}_0,\eqno(5.11)$$
$$\iota(\xi_i)=\ptl_{x_i},\;\;\iota(\eta_i)=P_i+\sum_{r=1}^{16}x_r\varrho([\xi_r,\eta_i])\qquad\for
\;\;i\in\ol{1,16}.\eqno(5.12)$$

In order to calculate $\{\iota(P_1),...,\iota(P_{16})\}$ explicitly,
we need the more formulas of $\nu$ on the positive root vectors of
$o(10,\mbb{C})$ extended from (2.42)-(2.47). We calculate
$$\nu(E_{3,5}-E_{10,8})=E_{\al_4+\al_5},\;\;\nu(E_{2,5}-E_{10,7})=-E_{\sum_{i=3}^5\al_i},\;\;\nu(E_{1,9}-E_{4,6})=E_{\sum_{i=1}^5\al_i},\eqno(5.13)$$
$$\nu(E_{3,9}-E_{4,8})=E_{\al_2+\al_4+\al_5},\;\;\nu(E_{2,9}-E_{4,7})=-E_{\sum_{i=2}^5\al_i},
\eqno(5.14)$$
$$\nu(E_{1,5}-E_{10,6})=E_{\al_1+\sum_{i=3}^5\al_i},\;\;\nu(E_{2,8}-E_{3,7})=-E_{\al_4+\sum_{i=2}^5\al_i},\eqno(5.15)$$
$$\nu(E_{1,8}-E_{3,6})=E_{\al_4+\sum_{i=1}^5\al_i},\;\;\nu(E_{1,7}-E_{2,6})
=-E_{\al_3+\al_4+\sum_{i=1}^5\al_i}.\eqno(5.16)$$

Now
\begin{eqnarray*} \iota(\eta_1)&=&P_1-x_1\varrho(\al_6)-x_2\varrho(E_{\al_5})
-x_3\varrho(E_{\al_4+\al_5})-x_4\varrho(E_{\al_3+\al_4+\al_5})
\\ &&-x_5\varrho(E_{\al_2+\al_4+\al_5})-x_6\varrho(E_{\sum_{i=2}^5\al_i})
-x_7\varrho(E_{\al_1+\sum_{i=3}^5\al_i})-x_8\varrho(E_{\sum_{i=1}^5\al_i})\\
&
&-x_9\varrho(E_{\al_4+\sum_{i=2}^5\al_i})-x_{10}\varrho(E_{\al_4+\sum_{i=1}^5\al_i})
-x_{12}\varrho(E_{\al_3+\al_4+\sum_{i=1}^5\al_i})\\ &=&
P_1+\frac{x_1}{2}[\sum_{i=1}^4(E_{i,i}-E_{5+i,5+i})-E_{5,5}+E_{10,10}-\kappa]
-x_2(E_{4,5}-E_{10,9})\\
&
&-x_3(E_{3,5}-E_{10,8})+x_4(E_{2,5}-E_{10,7})-x_5(E_{3,9}-E_{4,8})+x_6(E_{2,9}-E_{4,7})
\\ & &-x_7(E_{1,5}-E_{10,6})-x_8(E_{1,9}-E_{4,6})+x_9(E_{2,8}-E_{3,7})\\ & &
-x_{10}(E_{1,8}-E_{3,6})+x_{12}(E_{1,7}-E_{2,6})\hspace{6.7cm}(5.17)\end{eqnarray*}
by (2.43), (2.46), (2.47), (4.2), (5.10) and (5.13)-(5.16).
Moreover,
\begin{eqnarray*}
\iota(\eta_2)&=&\iota([E_{-\al_5},\eta_1])=[\iota(E_{-\al_5}),\iota(\eta_1)]
=-[(E_{5,4}-E_{9,10})|_{\cal A}+(E_{5,4}-E_{9,10}),\iota(\eta_1)]
\\&=&-[-x_2\ptl_{x_1}+x_{11}\ptl_{x_9}+x_{13}\ptl_{x_{10}}+x_{14}\ptl_{x_{12}}+(E_{5,4}-E_{9,10}),\iota(\eta_1)]
\\ &=&
P_2+\frac{x_2}{2}[\sum_{i\neq
4}(E_{i,i}-E_{5+i,5+i})-E_{4,4}+E_{9,9}-\kappa]-x_1(E_{5,4}-E_{9,10})
\\ &&
-x_3(E_{3,4}-E_{9,8}) +x_4(E_{2,4}-E_{9,7})
+x_5(E_{3,10}-E_{5,8})-x_6(E_{2,10}-E_{5,7})
\\ & &-x_7(E_{1,4}-E_{9,6})+x_8(E_{1,10}-E_{5,6})-x_{11}(E_{2,8}-E_{3,7})\\ & &
+x_{13}(E_{1,8}-E_{3,6})-x_{14}(E_{1,7}-E_{2,6})\hspace{6.7cm}(5.18)\end{eqnarray*}
by (3.12). According to (3.11),
\begin{eqnarray*}
\iota(\eta_3)&=&\iota([E_{-\al_4},\eta_2])=[\iota(E_{-\al_4}),\iota(\eta_2)]
=-[(E_{4,3}-E_{8,9})|_{\cal A}+(E_{4,3}-E_{8,9}),\iota(\eta_2)]
\\&=&-[-x_3\ptl_{x_2}-x_9\ptl_{x_6}-x_{10}\ptl_{x_8}+x_{15}\ptl_{x_{14}}+(E_{4,3}-E_{8,9}),\iota(\eta_2)]
\\ &=&
P_3+\frac{x_3}{2}[\sum_{i\neq
3}(E_{i,i}-E_{5+i,5+i})-E_{3,3}+E_{8,8}-\kappa]-x_1(E_{5,3}-E_{8,10})
\\ &&-x_2(E_{4,3}-E_{8,9})
 +x_4(E_{2,3}-E_{8,7}) -x_5(E_{4,10}-E_{5,9})-x_7(E_{1,3}-E_{8,6})\\ & &-x_9(E_{2,10}-E_{5,7})
+x_{10}(E_{1,10}-E_{5,6})+x_{11}(E_{2,9}-E_{4,7})\\ & &
-x_{13}(E_{1,9}-E_{4,6})+x_{15}(E_{1,7}-E_{2,6}).\hspace{6.6cm}(5.19)\end{eqnarray*}
Furthermore, (3.10) implies
\begin{eqnarray*}
\iota(\eta_4)&=&-\iota([E_{-\al_3},\eta_3])=-[\iota(E_{-\al_3}),\iota(\eta_3)]
=[(E_{3,2}-E_{7,8})|_{\cal A}+(E_{3,2}-E_{7,8}),\iota(\eta_3)]
\\&=&[x_4\ptl_{x_3}+x_6\ptl_{x_5}+x_{12}\ptl_{x_{10}}+x_{14}\ptl_{x_{13}}
+(E_{3,2}-E_{7,8}),\iota(\eta_3)]
\\ &=&
P_4+\frac{x_4}{2}[\sum_{i\neq
2}(E_{i,i}-E_{5+i,5+i})-E_{2,2}+E_{7,7}-\kappa]+x_1(E_{5,2}-E_{7,10})
\\ &&+x_2(E_{4,2}-E_{7,9})
 +x_3(E_{3,2}-E_{7,8}) -x_6(E_{4,10}-E_{5,9})+x_7(E_{1,2}-E_{7,6})\\ & &-x_9(E_{3,10}-E_{5,8})
+x_{11}(E_{3,9}-E_{4,8})+x_{12}(E_{1,10}-E_{5,6})\\ & &
-x_{14}(E_{1,9}-E_{4,6})+x_{15}(E_{1,8}-E_{3,6}).\hspace{6.6cm}(5.20)\end{eqnarray*}

Observe
\begin{eqnarray*} \iota(\eta_5)&=&\iota([E_{-\al_2},\eta_3])=[\iota(E_{-\al_2}),\iota(\eta_3)]
=-[(E_{10,4}-E_{9,5})|_{\cal A}+(E_{10,4}-E_{9,5}),\iota(\eta_3)]
\\&=&-[-x_5\ptl_{x_3}-x_6\ptl_{x_4}-x_8\ptl_{x_7}+x_{16}\ptl_{x_{15}}
+E_{10,4}-E_{9,5},\iota(\eta_3)]
\\ &=&
P_5+\frac{x_5}{2}[\sum_{i=1}^2(E_{i,i}-E_{5+i,5+i})-\sum_{i=
3}^5(E_{i,i}-E_{5+i,5+i})-\kappa]-x_1(E_{9,3}-E_{8,4})
\\ &&+x_2(E_{10,3}-E_{8,5})-x_3(E_{10,4}-E_{9,5})
 +x_6(E_{2,3}-E_{8,7}) -x_8(E_{1,3}-E_{8,6})\\ & &-x_9(E_{2,4}-E_{9,7})
+x_{10}(E_{1,4}-E_{9,6})-x_{11}(E_{2,5}-E_{10,7})\\ & &
+x_{13}(E_{1,5}-E_{10,6})-x_{16}(E_{1,7}-E_{2,6})\hspace{6.5cm}(5.21)\end{eqnarray*}
by (3.13). Similarly,
\begin{eqnarray*} \iota(\eta_6)&=&\iota([E_{-\al_2},\eta_4])=[\iota(E_{-\al_2}),\iota(\eta_3)]
=-[(E_{10,4}-E_{9,5})|_{\cal A}+(E_{10,4}-E_{9,5}),\iota(\eta_4)]
\\&=&-[-x_5\ptl_{x_3}-x_6\ptl_{x_4}-x_8\ptl_{x_7}+x_{16}\ptl_{x_{15}}
+E_{10,4}-E_{9,5},\iota(\eta_4)]
\\ &=&
P_6+\frac{x_6}{2}[\sum_{i=1,3}(E_{i,i}-E_{5+i,5+i})-\sum_{i=2,4,5
}(E_{i,i}-E_{5+i,5+i})-\kappa] +x_1(E_{9,2}-E_{7,4})
\\ &&-x_2(E_{10,2}-E_{7,5})
 -x_4(E_{10,4}-E_{9,5})+x_5(E_{3,2}-E_{7,8}) +x_8(E_{1,2}-E_{7,6})\\ & &-x_9(E_{3,4}-E_{9,8})
-x_{11}(E_{3,5}-E_{10,8})+x_{12}(E_{1,4}-E_{9,6})\\ & &
+x_{14}(E_{1,5}-E_{10,6})-x_{16}(E_{1,8}-E_{3,6}).\hspace{6.4cm}(5.22)\end{eqnarray*}
Moreover, (3.9) yields
\begin{eqnarray*}
\iota(\eta_7)&=&-\iota([E_{-\al_1},\eta_4])=-[\iota(E_{-\al_1}),\iota(\eta_4)]
=[(E_{2,1}-E_{6,7})|_{\cal A}+(E_{2,1}-E_{6,7}),\iota(\eta_4)]
\\&=&[x_7\ptl_{x_4}+x_8\ptl_{x_6}+x_{10}\ptl_{x_9}+x_{13}\ptl_{x_{11}}
+(E_{2,1}-E_{6,7}),\iota(\eta_4)]
\\ &=&
P_7+\frac{x_7}{2}[\sum_{i
2}^5(E_{i,i}-E_{5+i,5+i})-E_{1,1}+E_{6,6}-\kappa]-x_1(E_{5,1}-E_{6,10})
\\ &&-x_2(E_{4,1}-E_{6,9})
 -x_3(E_{3,1}-E_{6,8}) +x_4(E_{2,1}-E_{6,7})-x_8(E_{4,10}-E_{5,9})\\ & &-x_{10}(E_{3,10}-E_{5,8})
+x_{12}(E_{2,10}-E_{5,7})+x_{13}(E_{3,9}-E_{4,8})\\ & &
-x_{14}(E_{2,9}-E_{4,7})+x_{15}(E_{2,8}-E_{3,7}).\hspace{6.6cm}(5.23)\end{eqnarray*}
Furthermore, (3.13) gives
\begin{eqnarray*}
\iota(\eta_8)&=&\iota([E_{-\al_2},\eta_7])=[\iota(E_{-\al_2}),\iota(\eta_7)]
=-[(E_{10,4}-E_{9,5})|_{\cal A}+(E_{10,4}-E_{9,5}),\iota(\eta_7)]
\\&=&-[-x_5\ptl_{x_3}-x_6\ptl_{x_4}-x_8\ptl_{x_7}+x_{16}\ptl_{x_{15}}
+E_{10,4}-E_{9,5},\iota(\eta_7)]
\\ &=&
P_8+\frac{x_8}{2}[\sum_{i=2,3}(E_{i,i}-E_{5+i,5+i})-\sum_{i=1,4,5
}(E_{i,i}-E_{5+i,5+i})-\kappa]-x_1(E_{9,1}-E_{6,4})
\\ &&+x_2(E_{10,1}-E_{6,5})
 -x_5(E_{3,1}-E_{6,8}) +x_6(E_{2,1}-E_{6,7})
 -x_7(E_{10,4}-E_{9,5})\\ & &-x_{10}(E_{3,4}-E_{9,8})
+x_{12}(E_{2,4}-E_{9,7})-x_{13}(E_{3,5}-E_{10,8})\\ & &
+x_{14}(E_{2,5}-E_{10,7})-x_{16}(E_{2,8}-E_{3,7}).\hspace{6.4cm}(5.24)\end{eqnarray*}

According to (3.11),
\begin{eqnarray*}
\iota(\eta_9)&=&\iota([E_{-\al_4},\eta_6])=[\iota(E_{-\al_4}),\iota(\eta_6)]
=-[(E_{4,3}-E_{8,9})|_{\cal A}+(E_{4,3}-E_{8,9}),\iota(\eta_6)]
\\&=&-[-x_3\ptl_{x_2}-x_9\ptl_{x_6}-x_{10}\ptl_{x_8}+x_{15}\ptl_{x_{14}}+(E_{4,3}-E_{8,9}),\iota(\eta_6)]
\\ &=&
P_9+\frac{x_9}{2}[\sum_{i=1,4}(E_{i,i}-E_{5+i,5+i})-\sum_{i=2,3,5
}(E_{i,i}-E_{5+i,5+i})-\kappa] +x_1(E_{8,2}-E_{7,3})
\\ &&-x_3(E_{10,2}-E_{7,5})
 -x_4(E_{10,3}-E_{8,5})-x_5(E_{4,2}-E_{7,9})-x_6(E_{4,3}-E_{8,9}))\\ & &
 +x_{10}(E_{1,2}-E_{7,6})
+x_{11}(E_{4,5}-E_{10,9})+x_{12}(E_{1,3}-E_{8,6})\\ & &
-x_{15}(E_{1,5}-E_{10,6})+x_{16}(E_{1,9}-E_{4,6})\hspace{6.5cm}(5.25)\end{eqnarray*}
and
\begin{eqnarray*}
\iota(\eta_{10})&=&\iota([E_{-\al_4},\eta_8])=[\iota(E_{-\al_4}),\iota(\eta_8)]
=-[(E_{4,3}-E_{8,9})|_{\cal A}+(E_{4,3}-E_{8,9}),\iota(\eta_6)]
\\&=&-[-x_3\ptl_{x_2}-x_9\ptl_{x_6}-x_{10}\ptl_{x_8}+x_{15}\ptl_{x_{14}}+(E_{4,3}-E_{8,9}),\iota(\eta_8)]
\\ &=&
P_{10}+\frac{x_{10}}{2}[\sum_{i=2,4}(E_{i,i}-E_{5+i,5+i})-\sum_{i=1,3,5
}(E_{i,i}-E_{5+i,5+i})-\kappa]-x_1(E_{8,1}-E_{6,3})
\\ &&+x_3(E_{10,1}-E_{6,5})
 +x_5(E_{4,1}-E_{6,9})
 -x_7(E_{10,3}-E_{8,5})-x_8(E_{4,3}-E_{8,9})\\ & &+x_9(E_{2,1}-E_{6,7})
+x_{12}(E_{2,3}-E_{8,7})+x_{13}(E_{4,5}-E_{10,9})\\ & &
-x_{15}(E_{2,5}-E_{10,7})+x_{16}(E_{2,9}-E_{4,7}).\hspace{6.2cm}(5.26)\end{eqnarray*}
Moreover, (3.12) yields
\begin{eqnarray*}
\iota(\eta_{11})&=&-\iota([E_{-\al_5},\eta_9])=-[\iota(E_{-\al_5}),\iota(\eta_9)]
=-[(E_{5,4}-E_{9,10})|_{\cal A}+(E_{5,4}-E_{9,10}),\iota(\eta_1)]
\\&=&[-x_2\ptl_{x_1}+x_{11}\ptl_{x_9}+x_{13}\ptl_{x_{10}}+x_{14}\ptl_{x_{12}}+(E_{5,4}-E_{9,10}),\iota(\eta_9)]
\\
&=&P_{11}+\frac{x_{11}}{2}[\sum_{i=1,5}(E_{i,i}-E_{5+i,5+i})-\sum_{i=2,3,4
}(E_{i,i}-E_{5+i,5+i})-\kappa] -x_2(E_{8,2}-E_{7,3})
\\ &&+x_3(E_{9,2}-E_{7,4})
 +x_4(E_{9,3}-E_{8,4})-x_5(E_{5,2}-E_{7,10})-x_6(E_{5,3}-E_{8,10}))\\ & &
 +x_9(E_{5,4}-E_{9,10})+x_{13}(E_{1,2}-E_{7,6})
+x_{14}(E_{1,3}-E_{8,6})\\ & &
+x_{15}(E_{1,4}-E_{9,6})+x_{16}(E_{1,10}-E_{5,6})\hspace{6.3cm}(5.27)\end{eqnarray*}
Furthermore, (3.10) implies
\begin{eqnarray*}
\iota(\eta_{12})&=&-\iota([E_{-\al_3},\eta_{10}])=-[\iota(E_{-\al_3}),\iota(\eta_{10})]
=[(E_{3,2}-E_{7,8})|_{\cal A}+(E_{3,2}-E_{7,8}),\iota(\eta_3)]
\\&=&[x_4\ptl_{x_3}+x_6\ptl_{x_5}+x_{12}\ptl_{x_{10}}+x_{14}\ptl_{x_{13}}
+(E_{3,2}-E_{7,8}),\iota(\eta_{10})]
\\ &=&
P_{12}+\frac{x_{12}}{2}[\sum_{i=3,4}(E_{i,i}-E_{5+i,5+i})-\sum_{i=1,2,5
}(E_{i,i}-E_{5+i,5+i})-\kappa]+x_1(E_{7,1}-E_{6,2})
\\ &&+x_4(E_{10,1}-E_{6,5})
 +x_6(E_{4,1}-E_{6,9})
 +x_7(E_{10,2}-E_{7,5})+x_8(E_{4,2}-E_{7,9})\\ & &+x_9(E_{3,1}-E_{6,8})
+x_{10}(E_{3,2}-E_{7,8})+x_{14}(E_{4,5}-E_{10,9})\\ & &
-x_{15}(E_{3,5}-E_{10,8})+x_{16}(E_{3,9}-E_{4,8}).\hspace{6.2cm}(5.28)\end{eqnarray*}

Note that (3.9) gives
\begin{eqnarray*}
\iota(\eta_{13})&=&-\iota([E_{-\al_1},\eta_{11}])=-[\iota(E_{-\al_1}),\iota(\eta_{11})]
=[(E_{2,1}-E_{6,7})|_{\cal A}+(E_{2,1}-E_{6,7}),\iota(\eta_4)]
\\&=&[x_7\ptl_{x_4}+x_8\ptl_{x_6}+x_{10}\ptl_{x_9}+x_{13}\ptl_{x_{11}}
+(E_{2,1}-E_{6,7}),\iota(\eta_{11})]
\\ &=&
P_{13}+\frac{x_{13}}{2}[\sum_{i=2,5}(E_{i,i}-E_{5+i,5+i})-\sum_{i=1,3,4
}(E_{i,i}-E_{5+i,5+i})-\kappa] +x_2(E_{8,1}-E_{6,3})
\\ &&-x_3(E_{9,1}-E_{6,4})
 +x_5(E_{5,1}-E_{6,10})+x_7(E_{9,3}-E_{8,4})-x_8(E_{5,3}-E_{8,10}))\\ & &
 +x_{10}(E_{5,4}-E_{9,10})+x_{11}(E_{2,1}-E_{6,7})
+x_{14}(E_{2,3}-E_{8,7})\\ & &
+x_{15}(E_{2,4}-E_{9,7})+x_{16}(E_{2,10}-E_{5,7}).\hspace{6.2cm}(5.29)\end{eqnarray*}
Moreover, (3.10) yields
\begin{eqnarray*}
\iota(\eta_{14})&=&-\iota([E_{-\al_3},\eta_{13}])=-[\iota(E_{-\al_3}),\iota(\eta_{13})]
=[(E_{3,2}-E_{7,8})|_{\cal A}+(E_{3,2}-E_{7,8}),\iota(\eta_{13})]
\\&=&[x_4\ptl_{x_3}+x_6\ptl_{x_5}+x_{12}\ptl_{x_{10}}+x_{14}\ptl_{x_{13}}
+(E_{3,2}-E_{7,8}),\iota(\eta_{13})]
\\ &=&
P_{14}+\frac{x_{14}}{2}[\sum_{i=3,5}(E_{i,i}-E_{5+i,5+i})-\sum_{i=1,2,4
}(E_{i,i}-E_{5+i,5+i})-\kappa] -x_2(E_{7,1}-E_{6,2})
\\ &&-x_4(E_{9,1}-E_{6,4})
 +x_6(E_{5,1}-E_{6,10})-x_7(E_{9,2}-E_{7,4})+x_8(E_{5,2}-E_{7,10}))\\ & &
 +x_{11}(E_{3,1}-E_{6,8})+x_{12}(E_{5,4}-E_{9,10})
+x_{13}(E_{3,2}-E_{7,8})\\ & &
+x_{15}(E_{3,4}-E_{9,8})+x_{16}(E_{3,10}-E_{5,8}).\hspace{6.2cm}(5.30)\end{eqnarray*}
Furthermore, (3.11) implies
\begin{eqnarray*}
\iota(\eta_{15})&=&-\iota([E_{-\al_4},\eta_{14}])=-[\iota(E_{-\al_4}),\iota(\eta_{14})]
=[(E_{4,3}-E_{8,9})|_{\cal A}+(E_{4,3}-E_{8,9}),\iota(\eta_{14})]
\\&=&[-x_3\ptl_{x_2}-x_9\ptl_{x_6}-x_{10}\ptl_{x_8}+x_{15}\ptl_{x_{14}}+(E_{4,3}-E_{8,9}),\iota(\eta_{14})]
\\ &=&
P_{15}+\frac{x_{15}}{2}[\sum_{i=4,5}(E_{i,i}-E_{5+i,5+i})-\sum_{i=1,2,3
}(E_{i,i}-E_{5+i,5+i})-\kappa]+x_3(E_{7,1}-E_{6,2})
\\ &&+x_4(E_{8,1}-E_{6,3})
 +x_7(E_{8,2}-E_{7,3})-x_9(E_{5,1}-E_{6,10})-x_{10}(E_{5,2}-E_{7,10}))\\ & &
 +x_{11}(E_{4,1}-E_{6,9})-x_{12}(E_{5,3}-E_{8,10})
+x_{13}(E_{4,2}-E_{7,9})\\ & &
+x_{14}(E_{4,3}-E_{8,9})+x_{16}(E_{4,10}-E_{5,9}).\hspace{6.2cm}(5.31)\end{eqnarray*}
Finally,  (3.13) gives
\begin{eqnarray*}
\iota(\eta_{16})&=&-\iota([E_{-\al_2},\eta_{15}])=-[\iota(E_{-\al_2}),\iota(\eta_{15})]
=[(E_{10,4}-E_{9,5})|_{\cal A}+(E_{10,4}-E_{9,5}),\iota(\eta_{15})]
\\&=&[-x_5\ptl_{x_3}-x_6\ptl_{x_4}-x_8\ptl_{x_7}+x_{16}\ptl_{x_{15}}
+E_{10,4}-E_{9,5},\iota(\eta_{15})]
\\ &=&
P_{16}-\frac{x_{16}}{2}[\sum_{i=1
}^5(E_{i,i}-E_{5+i,5+i})+\kappa]-x_5(E_{7,1}-E_{6,2})
-x_6(E_{8,1}-E_{6,3})\\ &&
 -x_8(E_{8,2}-E_{7,3})+x_9(E_{9,1}-E_{6,4})+x_{10}(E_{9,2}-E_{7,4}))
 +x_{11}(E_{10,1}-E_{6,5})\\ & &+x_{12}(E_{9,3}-E_{8,4})
+x_{13}(E_{10,2}-E_{7,5}) +x_{14}(E_{10,3}-E_{8,5})\\
& &+x_{15}(E_{10,4}-E_{9,5}).\hspace{9.5cm}(5.32)\end{eqnarray*}

Recall ${\cal A}=\mbb{C}[x_1,...,x_{16}]$. Let $M$ be an
$o(10,\mbb{C})$-module and set
$$\widehat{M}={\cal A}\otimes_{\mbb{C}}M.\eqno(5.33)$$
We identify
$$f\otimes v=fv\qquad\for\;\;f\in{\cal A},\;v\in M.\eqno(5.34)$$
Recall the Lie algebra $\widehat{\cal K}$ defined via (5.1)-(5.8).
Fix $c\in\mbb{C}$. Then $\widehat M$ becomes a $\widehat{\cal
K}$-module with the action defined by
$$d(fv)=d(f)v,\;\;\kappa|(fv)=cfv,\;\;(gB)(fv)=fg
B(v)\eqno(5.35)$$ for $f,g\in{\cal A},\;v\in M$ and $B\in
o(10,\mbb{C})$.

Since the linear map $\iota: {\cal G}^{E_6}\rta \widehat{\cal K}$
defined in (5.10)-(5.12) is a  Lie algebra monomorphism,
$\widehat{M}$ becomes a ${\cal G}^{E_6}$-module with the action
defined by
$$\xi(w)=\iota(\xi)(w)\qquad\for\;\;\xi\in{\cal
G}^{E_6},\;w\in\widehat{M}.\eqno(5.36)$$ In fact, we have:\psp

{\bf Theorem 5.1}. {\it The map $M\mapsto \widehat M$ gives a
functor from $o(10,\mbb{C})$-{\bf Mod} to ${\cal G}^{E_6}$-{\bf
Mod}.}\psp

We remark that  the module $\widehat M$ is not a generalized module
in general because it may not be equal to $U({\cal G})(M)=U({\cal
G}_-)(M)$. \psp

{\bf Proposition 5.2}. {\it If $M$ is an irreducible
$o(10,\mbb{C})$-module, then $U({\cal G}_-)(M)$ is an irreducible
${\cal G}^{E_6}$-module.}

{\it Proof.} Note that for any $i\in\ol{1,16}$, $f\in{\cal A}$ and
$v\in M$, (5.12), (5.35) and (5.36) imply
$$\xi_i(fv)=\ptl_{x_i}(f)v.\eqno(5.37)$$
Let $W$ be any nonzero ${\cal G}^{E_6}$-submodule. The above
expression shows that $W\bigcap M\neq\{0\}$. According to (5.35),
$W\bigcap M$ is an $o(10,\mbb{C})$-submodule. By the irreducibility
of $M$, $M\subset W$. Thus $U({\cal G}_-)(M)\subset W$. So $U({\cal
G}_-)(M)=W$ is irreducible. $\qquad\Box$\psp

By the above proposition, the map $M\mapsto U({\cal G}_-)(M)$ is a
polynomial extension from  irreducible $o(10,\mbb{C})$-modules to
irreducible ${\cal G}^{E_6}$-modules.

\section{Irreducibility}

In this section, we want to determine the irreducibility of ${\cal
G}^{E_6}$-modules.

   Let $M$ be an $o(10,\mbb{C})$-module and let $\widehat M$ be the
${\cal G}^{E_6}$-module defined in (5.33)-(5.36).   Moreover,
$\widehat{M}$ can be viewed as an $o(10,\mbb{C})$-module with the
representation $\iota(\nu(B))|_{\widehat{M}}$ (cf. (2.42)-(2.47)).
Indeed, (5.11) and (5.36) show
$$\nu(B)(fv)=B(f)v+fB(v)\qquad\for\;\;B\in o(10,\mbb{C}),\;f\in{\cal
A},\;v\in M\eqno(6.1)$$ (cf. (2.49)-(2.71)). So $\widehat{M}={\cal
A}\otimes_{\mbb{C}}M$ is a tensor module of $o(10,\mbb{C})$. Write
$$\eta^\al=\prod_{i=1}^{16}\eta_{i}^{\alpha_i},\;\;
 |\al|=\sum_{i=1}^{16}\al_i\;\; \for\;\;
\al=(\al_1,\al_2,...,\al_{16})\in\mbb{N}^{16}\eqno(6.2)$$ (cf.
(2.27)-(2.31)).
 Recall the Lie subalgebras ${\cal G}_\pm$ and ${\cal G}_0$ of ${\cal G}^{E_6}$ defined
 in (2.32).  For $k\in\mbb{N}$, we set
$${\cal A}_k=\mbox{Span}_{\mathbb{C}}\{x^\al\mid
\al \in\mbb{N}^{16};|\al|=k\},\;\; \widehat
M_{{\langle}k\rangle}={\cal A}_kM\eqno(6.3)$$ (cf. (2.39), (5.34))
and
 $$(U({\cal G}_-)(
M))_{{\langle}k\rangle}=\mbox{Span}_{\mathbb{C}}\{\eta^\al( M)\mid
\al \in\mbb{N}^{16}, \ |\al|=k\}. \eqno(6.4)$$
 Moreover,
$$(U({\cal G}_-)(
M))_{{\langle}0\rangle}=\widehat M_{{\langle}0\rangle}=
M.\eqno(6.5)$$ Furthermore,
 $$\widehat M=\bigoplus\limits_{k=0}^\infty\widehat M_{\langle
 k\rangle},\qquad
 U({\cal G}_-)(M)=\bigoplus\limits_{k=0}^\infty(U({\cal G}_-)(M))_{\langle k\rangle}.\eqno(6.6)$$

Next we define a linear transformation $\vf$ on  $\widehat M$
determined by
$$\vf(x^\al v)=\eta^\al(
v)\qquad\for\;\;\al\in\mbb{N}^{16},\;v\in M.\eqno(6.7)$$ Note that
${\cal A}_1=\sum_{i=1}^{16}\mbb{C}x_i$ forms the $16$-dimensional
${\cal G}_0$-module (equivalently, the $o(10,\mbb{C})$ spin module).
According to (2.9) and (2.10), ${\cal G}_-$ forms a ${\cal
G}_0$-module with respect to the adjoint representation, and the
linear map from ${\cal A}_1$ to ${\cal G}_0$ determined by
$x_i\mapsto \eta_i$ for $i\in\ol{1,16}$ gives a ${\cal G}_0$-module
isomorphism. Thus $\vf$ can also be viewed as a ${\cal G}_0$-module
homomorphism from $\widehat M$ to $U({\cal G}_0)(M)$. Moreover,
$$\vf(\widehat M_{\la k\ra})=(U({\cal G}_-)(
M))_{{\langle}k\rangle}\qquad\for\;\;k\in\mbb{N}.\eqno(6.8)$$

Note that the Casimir element of $o(10,\mbb{C})$ is
\begin{eqnarray*}\omega&=&\sum_{1\leq
i<j\leq
10}[(E_{i,5+j}-E_{j,5+i})(E_{5+j,i}-E_{5+i,j})+(E_{5+j,i}-E_{5+i,j})(E_{i,5+j}-E_{j,5+i})]
\\ & &+\sum_{i,j=1}^5(E_{i,j}-E_{5+j,5+i})(E_{j,i}-E_{5+i,5+j})\in
U(o(10,\mbb{C})).\hspace{4.2cm}(6.9)\end{eqnarray*}  The algebra
$U(o(10,\mbb{C}))$  can be imbedded into the tensor algebra
$U(o(10,\mbb{C}))\otimes U(o(10,\mbb{C}))$ by the associative
algebra homomorphism $\dlt: U(o(10,\mbb{C})) \rightarrow
U(o(10,\mbb{C}))\otimes_\mbb{C} U(o(10,\mbb{C}))$ determined  by
$$\dlt(u)=u\otimes 1 +1 \otimes u \qquad \mbox{ for} \ u\in
o(10,\mbb{C}).\eqno(6.10)$$ Set
$$\td\omega=\frac{1}{2}(\dlt(\omega)-\omega\otimes 1-1\otimes
\omega)\in U(o(10,\mbb{C}))\otimes_\mbb{C}
U(o(2n,\mbb{C})).\eqno(6.11)$$ By (6.9),
\begin{eqnarray*}\td\omega\!&=&\!\sum_{1\leq
i<j\leq
5}[(E_{i,5+j}-E_{j,5+i})\otimes(E_{5+j,i}-E_{5+i,j})+(E_{5+j,i}-E_{5+i,j})\otimes(E_{i,5+j}-E_{j,5+i})]
\\ & &+\sum_{i,j=1}^5(E_{i,j}-E_{5+j,5+i})\otimes(E_{j,i}-E_{5+i,5+j}).\hspace{6.2cm}(6.12)\end{eqnarray*}
Furthermore, $\td{\omega}$ acts on $\widehat M$ as an
$o(10,\mbb{C})$-module homomorphism via
$$(B_1\otimes B_2)(fv)=B_1(f)B_2(v)\qquad\for\;\;B_1,B_2\in
o(10,\mbb{C}).\eqno(6.13)$$

{\bf Lemma 6.1}. {\it We have $\vf|_{\widehat M_{\la
1\ra}}=(\td\omega-c/2)|_{\widehat M_{\la 1\ra}}$.}

{\it Proof.} For any $v\in M$,  (2.49)-(2.71), (5.17), (6.12) and
(6.13) give
\begin{eqnarray*}\td\omega(x_1v)&=&[-x_2(E_{4,5}-E_{10,9})-x_3(E_{3,5}-E_{10,8})
+x_4(E_{2,5}-E_{10,7})\\&&-x_5(E_{3,9}-E_{4,8})-x_7(E_{1,5}-E_{10,6})
+x_6(E_{2,9}-E_{4,7})\\ &
&-x_8(E_{1,9}-E_{4,6})+x_9(E_{2,8}-E_{3,7})-x_{10}(E_{1,8}-E_{3,6})\\
& &
+x_{12}(E_{1,7}-E_{2,6})+x_1(\sum_{i=1}^4(E_{i,i}-E_{5+i,5+i})-E_{5,5}+E_{10,10})]v
\\ &=& \eta_1(v)+(c/2)x_1v=(\vf+c/2)(x_1v),\hspace{6.4cm}(6.14)
\end{eqnarray*}
equivalently, $\vf(x_1v)=(\td\omega-c/2)(x_1v)$. According to
(3.12), $(E_{9,10}-E_{5,4})(x_1)=x_2$. So
\begin{eqnarray*}\qquad& &\vf(x_2v)+\vf(x_1(E_{9,10}-E_{5,4})(v))=
(E_{9,10}-E_{5,4})(\vf(x_1v)\\
&=&(E_{9,10}-E_{5,4})[(\td\omega-c/2)(x_1v)]=(\td\omega-c/2)[(E_{9,10}-E_{5,4})(x_1v)]
\\
&=&(\td\omega-c/2)[x_2v+x_1(E_{9,10}-E_{5,4})(v)]\\
&=&(\td\omega-c/2)(x_2v)+(\td\omega-c/2)(x_1(E_{9,10}-E_{5,4})(v))
\\ &=&(\td\omega-c/2)(x_2v)+\vf(x_1(E_{9,10}-E_{5,4})(v)),\hspace{5.8cm}(6.15)\end{eqnarray*}
equivalently, $\vf(x_2v)=(\td\omega-c/2)(x_2v)$.

Observe that
$$(E_{8,9}-E_{4,3})(x_2)=x_3,\;\;(E_{3,2}-E_{7,8})(x_3)=x_4,\;\;(E_{9,5}-E_{10,4})(x_3)=x_5,\eqno(6.16)$$
$$(E_{6,8}-E_{3,1})(x_3)=x_7,\;\;(E_{7,5}-E_{10,2})(x_2)=x_6,\;\;(E_{10,1}-E_{6,5})(x_2)=x_8,\eqno(6.17)$$
$$(E_{8,2}-E_{7,3})(x_1)=x_9,\;\;(E_{6,3}-E_{8,1})(x_1)=x_{10},\;\;(E_{7,3}-E_{8,2})(x_2)=x_{11},\eqno(6.17)$$
$$(E_{3,2}-E_{7,8})(x_{10})=x_{12},\;\;(E_{8,1}-E_{6,3})(x_2)=x_{13},\;\;(E_{3,2}-E_{7,8})(x_{13})=x_{14},\eqno(6.18)$$
$$(E_{4,3}-E_{8,9})(x_{14})=x_{15},\;\;(E_{10,4}-E_{9,5})(x_{15})=x_{16}\eqno(6.19)$$
by (3.9)-(3.13). Using the argument similar to (6.14) and induction,
we can prove
$$\vf(x_iv)=(\td\omega-c/2)(x_iv)\qquad\for\;\;i\in\ol{1,16},\eqno(6.20)$$
equivalently, the lemma holds. $\qquad\Box$\psp

We take the subspace
$${\cal H}=\sum_{i=1}^5\mbb{C}(E_{i,i}-E_{n+i,n+i})\eqno(6.21)$$
as a Cartan subalgebra of the Lie algebra $o(10,\mbb{C})$  and
define $\{\ves_1,...,\ves_5\}\subset{\cal H}^\ast$ by
$$\ves_i(E_{j,j}-E_{n+j,n+j})=\dlt_{i,j}.\eqno(6.22)$$
The inner product $(\cdot,\cdot)$ on the $\mbb{Q}$-subspace
$$L_\mbb{Q}=\sum_{i=1}^5\mbb{Q}\ves_i\eqno(6.23)$$
is given by
$$(\ves_i,\ves_j)=\dlt_{i,j}\qquad\for\;\;i,j\in\ol{1,n}.\eqno(6.24)$$
Then the root system of $o(10,\mbb{C})$ is
$$\Phi_{D_5}=\{\pm \ves_i\pm\ves_j\mid1\leq i<j\leq
5\}.\eqno(6.25)$$ We take the set of positive roots
$$\Phi_{D_5}^+=\{\ves_i\pm\ves_j\mid1\leq i<j\leq
5\}.\eqno(6.26)$$ In particular,
$$\Pi_{D_5}=\{\ves_1-\ves_2,...,\ves_4-\ves_5,\ves_4+\ves_5\}\;\;\mbox{is the set of positive simple roots}.\eqno(6.27)$$

Recall the set of dominate integral weights
$$\Lmd^+=\{\mu\in L_\mbb{Q}\mid
(\ves_4+\ves_5,\mu),(\ves_i-\ves_{i+1},\mu)\in\mbb{N}\;\for\;i\in\ol{1,4}\}.\eqno(6.28)$$
According to (6.24),
$$\Lmd^+=\{\mu=\sum_{i=1}^5\mu_i\ves_i\mid
\mu_i\in\mbb{Z}/2;\mu_i-\mu_{i+1},\mu_4+\mu_5\in\mbb{N}\}.\eqno(6.29)$$
Note that if $\mu\in\Lmd^+$, then $\mu_4\geq |\mu_5|$. Denote
$$\rho=\frac{1}{2}\sum_{\nu\in\Phi_{D_5}^+}\nu.\eqno(6.30)$$
Then
$$(\rho,\nu)=1\qquad\for\;\;\nu\in\Pi_{D_5}\eqno(6.31)$$
(e.g., cf. [Hu]). By (6.24),
$$\rho=\sum_{i=1}^4(5-i)\ves_i.\eqno(6.32)$$
For any $\mu\in \Lmd^+$, we denote by $V(\mu)$ the
finite-dimensional irreducible $o(10,\mbb{C})$-module with the
highest weight $\mu$ and have
$$\omega|_{V(\mu)}=(\mu+2\rho,\mu)\mbox{Id}_{V(\mu)}\eqno(6.33)$$
by (6.9).

Let $\mbb{Z}_2=\mbb{Z}/2\mbb{Z}=\{0,1\}$. According to (2.71) and
Table 1, the weight set of the $o(10,\mbb{C})$-module ${\cal A}_1$
is
$$\Pi({\cal A}_1)=\{(1/2)\sum_{i=1}^5(-1)^{k_i}\ves_i\mid
k_i\in\mbb{Z}_2,\;\sum_{i=1}^5k_i=1\}.\eqno(6.34)$$ Fix
$\lmd\in\Lmd^+$, we define
$$\Upsilon(\lmd)=\{\lmd+\mu\mid\mu\in\Pi({\cal
A}_1),\;\lmd+\mu\in\Lmd^+\}.\eqno(6.35)$$ \pse

{\bf Lemma 6.2}. {\it We have}:
$${\cal A}_1\otimes V(\lmd)\cong \bigoplus_{\lmd'\in
\Upsilon(\lmd)}V(\lmd').\eqno(6.36)$$

{\it Proof}. Note that all the weight subspaces of ${\cal A}_1$ are
one-dimensional. Thus all the irreducible components of ${\cal
A}_1\otimes V(\lmd)$ are of multiplicity one. Since
$$\rho+\lmd+\mu\in\Lmd^+\qquad\for\;\;\mu\in \Pi({\cal
A}_1),\eqno(6.37)$$ the tensor theory of finite-dimensional
irreducible modules over a finite-dimensional simple Lie algebra
(e.g., cf. [Hu]) says that $V(\lmd')$ is a component of ${\cal
A}_1\otimes V(\lmd)$ if and only if $\lmd'\in
\Upsilon(\lmd).\qquad\Box$\psp

Recall
$$\mbox{the highest weight of}\;{\cal
A}_1=\frac{1}{2}(\ves_1+\ves_2+\ves_3+\ves_4-\ves_5)=\lmd_4,\eqno(6.38)$$
the forth fundamental weight of $o(10,\mbb{C})$, by (2.71) and Table
1. Thus the eigenvalues of
$\td{\omega}|_{\widehat{V(\lmd)}_{\la1\ra}}$ are
$$\{[(\lmd'+2\rho,\lmd')-(\lmd+2\rho,\lmd)-(\lmd_4+2\rho,\lmd_4)]/2\mid\lmd'\in
\Upsilon(\lmd)\}\eqno(6.39)$$ by (6.11) and (6.13). Define
$$\ell_\omega(\lmd)=\min\{[(\lmd'+2\rho,\lmd')-(\lmd+2\rho,\lmd)-(\lmd_4+2\rho,\lmd_4)]/2\mid\lmd'\in
\Upsilon(\lmd)\},\eqno(6.40)$$ which will be used to determine the
irreducibility of $\widehat{V(\lmd)}$. If
$\lmd'=\lmd+\lmd_4-\al\in\Upsilon(\lmd)$ with $\al\in\Phi_{D_5}^+$,
then
$$(\lmd'+2\rho,\lmd')-(\lmd+2\rho,\lmd)-(\lmd_4+2\rho,\lmd_4)
=2[(\lmd,\lmd_4)+1-(\rho+\lmd+\lmd_4,\al)].\eqno(6.41)$$

Recall the differential operators $P_1,...,P_{16}$ given in
(4.16)-(4.32). We also view the elements of ${\cal A}$ as the
multiplication operators on ${\cal A}$. Recall $\zeta_1$ in (3.8).
It turns out that we need the following lemma in order to determine
the irreducibility of $\widehat{V(\lmd)}$.\psp

{\bf Lemma 6.3}. {\it As operators on} ${\cal A}$:
$$P_{11}x_1+P_1x_{11}+P_9x_2+P_2x_9-P_6x_3-P_3x_6+P_5x_4+P_4x_5=\zeta_1(D-6).\eqno(6.42)$$

{\it Proof}. According to  (4.16), (4.18)-(4.22), (4.25) and (4.27),
we find that
\begin{eqnarray*}\qquad&
&P_{11}x_1+P_1x_{11}+P_9x_2+P_2x_9-P_6x_3-P_3x_6+P_5x_4+P_4x_5\\
&=&-6\zeta_1+x_1P_{11}+x_{11}P_1+x_2P_9+x_9P_2-x_3P_6-x_6P_3+x_4P_5+x_5P_4\hspace{1.2cm}(6.43)\end{eqnarray*}
and
\begin{eqnarray*}&
&x_1P_{11}+x_{11}P_1+x_2P_9+x_9P_2-x_3P_6-x_6P_3+x_4P_5+x_5P_4 \\&=&
x_1(x_{11}D-\zeta_1\ptl_{x_1}+\zeta_5\ptl_{x_7}+\zeta_9\ptl_{x_8}+\zeta_8\ptl_{x_{10}}-\zeta_7\ptl_{x_{12}})
\\&
&+x_{11}(x_1D-\zeta_1\ptl_{x_{11}}-\zeta_2\ptl_{x_{13}}-\zeta_3\ptl_{x_{14}}-\zeta_4\ptl_{x_{15}}-\zeta_{10}\ptl_{x_{16}})
\\&&+x_2(x_9D-\zeta_1\ptl_{x_2}+\zeta_4\ptl_{x_7}-\zeta_{10}\ptl_{x_8}-\zeta_8\ptl_{x_{13}}+\zeta_7\ptl_{x_{14}})
\\&&+x_9(x_2D-\zeta_1\ptl_{x_9}-\zeta_2\ptl_{x_{10}}-\zeta_3\ptl_{x_{12}}+\zeta_5\ptl_{x_{15}}-\zeta_9\ptl_{x_{16}})
\\&&-x_3(x_6D+\zeta_1\ptl_{x_3}-\zeta_3\ptl_{x_7}+\zeta_{10}\ptl_{x_{10}}-\zeta_9\ptl_{x_{13}}+\zeta_7\ptl_{x_{15}})
\\&&-x_6(x_3D+\zeta_1\ptl_{x_6}+\zeta_2\ptl_{x_8}+\zeta_4\ptl_{x_{12}}+\zeta_5\ptl_{x_{14}}-\zeta_8\ptl_{x_{16}})\hspace{5cm}\end{eqnarray*}
\begin{eqnarray*}
&&+x_4(x_5D-\zeta_1\ptl_{x_4}-\zeta_2\ptl_{x_7}-\zeta_{10}\ptl_{x_{12}}-\zeta_9\ptl_{x_{14}}-\zeta_8\ptl_{x_{15}})
\\&&+x_5(x_4D-\zeta_1\ptl_{x_5}+\zeta_3\ptl_{x_8}-\zeta_4\ptl_{x_{10}}-\zeta_5\ptl_{x_{13}}+\zeta_7\ptl_{x_{16}})
\\&=&2\zeta_1D-\zeta_1\sum_{i=1,2,3,4,5,6,9,11}x_i\ptl_{x_i}+(x_1\zeta_5+x_2\zeta_4+x_3\zeta_3-x_4\zeta_2)\ptl_{x_7}
\\&&+(x_1\zeta_9-x_2\zeta_{10}-x_6\zeta_2+x_5\zeta_3)\ptl_{x_8}
+(x_1\zeta_8-x_9\zeta_2-x_3\zeta_{10}-x_5\zeta_4)\ptl_{x_{10}}
\\&&-(x_1\zeta_7+x_9\zeta_3+x_6\zeta_4+x_4\zeta_{10})\ptl_{x_{12}}
-(x_{11}\zeta_2+x_2\zeta_8-x_3\zeta_9+x_5\zeta_5)\ptl_{x_{13}}
\\&&-(x_{11}\zeta_3-x_2\zeta_7+x_6\zeta_5-x_4\zeta_9)\ptl_{x_{14}}
-(x_{11}\zeta_4-x_9\zeta_5+x_3\zeta_7+x_4\zeta_8)\ptl_{x_{15}}
\\&&-(x_{11}\zeta_{10}+x_9\zeta_9-x_6\zeta_8-x_5\zeta_7)\ptl_{x_{16}}
=\zeta_1D \hspace{6.5cm}(6.44)\end{eqnarray*}
 because\begin{eqnarray*}\qquad&
&x_1\zeta_5+x_2\zeta_4+x_3\zeta_3-x_4\zeta_2
=x_1(-x_2x_{15}-x_3x_{14}+x_4x_{13}-x_7x_{11})\\
& &+x_2(x_1x_{15}-x_3x_{12}+x_4x_{10}-x_7x_9)
+x_3(x_1x_{14}+x_2x_{12}-x_4x_8+x_6x_7)\\
& &-x_4(x_1x_{13}+x_2x_{10}-x_3x_8+x_5x_7)=-\zeta_1x_7
\hspace{5.7cm}(6.45)\end{eqnarray*} by (3.8) and (3.14)-(3.17),
\begin{eqnarray*}\qquad&
&x_1\zeta_9-x_2\zeta_{10}-x_6\zeta_2+x_5\zeta_3=x_1(x_2x_{16}-x_5x_{14}+x_6x_{13}-x_8x_{11})
\\&
&-x_2(x_1x_{16}+x_5x_{12}-x_6x_{10}+x_8x_9)-x_6(x_1x_{13}+x_2x_{10}-x_3x_8+x_5x_7)
\\& &+x_5(x_1x_{14}+x_2x_{12}-x_4x_8+x_6x_7)=-\zeta_1x_8\hspace{5.7cm}(6.46)\end{eqnarray*}
by (3.8), (3.14), (3.15), (3.18) and (3.19),
\begin{eqnarray*}\qquad&
&x_1\zeta_8-x_9\zeta_2-x_3\zeta_{10}+x_5\zeta_4=x_1(x_3x_{16}+x_5x_{15}+x_9x_{13}-x_{10}x_{11})\\
&
&-x_9(x_1x_{13}+x_2x_{10}-x_3x_8+x_5x_7)-x_3(x_1x_{16}+x_5x_{12}-x_6x_{10}+x_8x_9)
\\ & &-x_5(x_1x_{15}-x_3x_{12}-x_4x_{10}-x_7x_9)=-\zeta_1x_{10}\hspace{5.4cm}(6.47)\end{eqnarray*}
by (3.8), (3.14), (3.16), (3.18) and (3.20),
\begin{eqnarray*}\qquad&&
x_1\zeta_7+x_9\zeta_3+x_6\zeta_4+x_4\zeta_{10}=x_1(-x_4x_{16}-x_6x_{15}-x_9x_{14}+x_{11}x_{12})
\\&
&+x_9(x_1x_{14}+x_2x_{12}-x_4x_8+x_6x_7)+x_6(x_1x_{15}-x_3x_{12}+x_4x_{10}-x_7x_9)\\
&
&+x_4(x_1x_{16}+x_5x_{12}-x_6x_{10}+x_8x_9)=\zeta_1x_{12}\hspace{5.7cm}(6.48)\end{eqnarray*}
by (3.8), (3.15), (3.16), (3.18) and (3.21),
\begin{eqnarray*}\qquad&&x_{11}\zeta_2+x_2\zeta_8-x_3\zeta_9+x_5\zeta_5=
x_{11}(x_1x_{13}+x_2x_{10}-x_3x_8+x_5x_7)\\ &
&+x_2(x_3x_{16}+x_5x_{15}+x_9x_{13}-x_{10}x_{11})-x_3(x_2x_{16}-x_5x_{14}+x_6x_{13}-x_8x_{11})
\\& &+x_5(-x_2x_{15}-x_3x_{14}+x_4x_{13}-x_7x_{11})=\zeta_1x_{13}\hspace{5.2cm}(6.49)\end{eqnarray*}
by (3.8), (3.14), (3.17), (3.19) and (3.20),
\begin{eqnarray*}\qquad&&x_{11}\zeta_3-x_2\zeta_7+x_6\zeta_5-x_4\zeta_9=x_{11}(x_1x_{14}+x_2x_{12}-x_4x_8+x_6x_7)
\\ &
&-x_2(-x_4x_{16}-x_6x_{15}-x_9x_{14}+x_{11}x_{12})+x_6(-x_2x_{15}-x_3x_{14}+x_4x_{13}-x_7x_{11})
\\& &-x_4(x_2x_{16}-x_5x_{14}+x_6x_{13}-x_8x_{11})=\zeta_1x_{14}\hspace{5.6cm}(6.50)\end{eqnarray*}
by (3.8), (3.15), (3.17), (3.19) and (3.21),
\begin{eqnarray*}\qquad&&x_{11}\zeta_4-x_9\zeta_5+x_3\zeta_7+x_4\zeta_8=
x_{11}(x_1x_{15}-x_3x_{12}+x_4x_{10}-x_7x_9)\\ &
&-x_9(-x_2x_{15}-x_3x_{14}+x_4x_{13}-x_7x_{11})+x_3(-x_4x_{16}-x_6x_{15}-x_9x_{14}+x_{11}x_{12})
\\& &+x_4(x_3x_{16}+x_5x_{15}+x_9x_{13}-x_{10}x_{11})=\zeta_1x_{15}\hspace{5.5cm}(6.51)\end{eqnarray*}
by (3.8), (3.16), (3.17), (3.20) and (3.21),
\begin{eqnarray*}\qquad&&x_{11}\zeta_{10}+x_9\zeta_9-x_6\zeta_8-x_5\zeta_7=x_{11}(x_1x_{16}+x_5x_{12}-x_6x_{10}+x_8x_9)
\\
&&+x_9(x_2x_{16}-x_5x_{14}+x_6x_{13}-x_8x_{11})-x_6(x_3x_{16}+x_5x_{15}+x_9x_{13}-x_{10}x_{11})\\&
&-x_5(-x_4x_{16}-x_6x_{15}-x_9x_{14}+x_{11}x_{12})=\zeta_1x_{16}\hspace{5.1cm}(6.52)\end{eqnarray*}
by (3.8) and (3.18)-(3.21).$\qquad\Box$\psp

We define the multiplication
$$f(gv)=(fg)v\qquad\for\;\;f,g\in{\cal A},\;v\in M.\eqno(6.53)$$
Then (5.17)-(5.22), (5.25) and (5.27) gives
\begin{eqnarray*}\qquad&&
\sum_{r=1}^{16}x_r\varrho([\xi_r,x_1\eta_{11}+x_{11}\eta_1+x_2\eta_9+x_9\eta_2-x_3\eta_6-x_6\eta_3+x_4\eta_5+x_5\eta_4])
\\ &=& x_1[x_{13}(E_{1,2}-E_{7,6}) +x_{14}(E_{1,3}-E_{8,6})
+x_{15}(E_{1,4}-E_{9,6})+x_{16}(E_{1,10}-E_{5,6})] \\ & &
-x_{11}[x_7(E_{1,5}-E_{10,6})+x_8(E_{1,9}-E_{4,6})+x_{10}(E_{1,8}-E_{3,6})-x_{12}(E_{1,7}-E_{2,6})]
\\ & &+x_2[x_{10}(E_{1,2}-E_{7,6})+x_{12}(E_{1,3}-E_{8,6})
-x_{15}(E_{1,5}-E_{10,6})+x_{16}(E_{1,9}-E_{4,6})]
\\ &&-x_9[x_7(E_{1,4}-E_{9,6})-x_8(E_{1,10}-E_{5,6})
-x_{13}(E_{1,8}-E_{3,6})+x_{14}(E_{1,7}-E_{2,6})]
\\& &-x_3[ x_8(E_{1,2}-E_{7,6})+x_{12}(E_{1,4}-E_{9,6})
+x_{14}(E_{1,5}-E_{10,6})-x_{16}(E_{1,8}-E_{3,6})]
\\ &&+x_6[x_7(E_{1,3}-E_{8,6})
-x_{10}(E_{1,10}-E_{5,6})+x_{13}(E_{1,9}-E_{4,6})-x_{15}(E_{1,7}-E_{2,6})]
\\& &-x_4[x_8(E_{1,3}-E_{8,6})
-x_{10}(E_{1,4}-E_{9,6})
-x_{13}(E_{1,5}-E_{10,6})+x_{16}(E_{1,7}-E_{2,6})]
\\& &+x_5[x_7(E_{1,2}-E_{7,6})+x_{12}(E_{1,10}-E_{5,6})
-x_{14}(E_{1,9}-E_{4,6})+x_{15}(E_{1,8}-E_{3,6})]
\\& &+\zeta_1(E_{1,1}-E_{6,6}-c)
\\&=&\sum_{i=1}^5\zeta_i(E_{1,i}-E_{5+i,6})+\sum_{r=2}^5\zeta_{5+r}(E_{1,5+r}-E_{r,6})-c\zeta_1
\hspace{3.9cm}(6.54)\end{eqnarray*} as operators on $\widehat M$
(cf. (5.33)), where $\zeta_i$ are defined in (3.8) and
(3.14)-(3.22).
 By Lemma 6.3, (5.5), (5.6) and (6.54),
\begin{eqnarray*}T_1&=&\iota(\eta_{11})x_1+\iota(\eta_1)x_{11}+\iota(\eta_9)x_2+\iota(\eta_2)x_9-\iota(\eta_6)x_3
-\iota(\eta_3)x_6+\iota(\eta_5)x_4+\iota(\eta_4)x_5
\\
&=&\zeta_1(D-c-6)+\sum_{i=1}^5\zeta_i(E_{1,i}-E_{5+i,6})+\sum_{r=2}^5\zeta_{5+r}(E_{1,5+r}-E_{r,6})
\hspace{2.5cm}(6.55)\end{eqnarray*}
 as operators on $\widehat M$. We define an $o(10,\mbb{C})$-module
 structure on the space $\mbox{End}\:\widehat M$ of linear
 transformations on $\widehat M$ by
 $$B(T)=[\nu(B),T]=\nu(B)T-T\nu(B)\qquad\for\;\;B\in
 o(10,\mbb{C}),\;T\in \mbox{End}\:\widehat M\eqno(6.56)$$
(cf. (6.1)). It can be verified that $T_1$ is an
$o(10,\mbb{C})$-singular vector with weight $\ves_1$ in
$\mbox{End}\:\widehat M$. So it generates the 10-dimensional natural
module.

  We set
\begin{eqnarray*}T_2&=&-[\iota(E_{-\al_1}),T_1]=[(E_{2,1}-E_{6,7})|_{\cal
A}+(E_{2,1}-E_{6,7}),T_1]
\\&=&\zeta_2(D-c-6)+\sum_{i=1}^5\zeta_i(E_{2,i}-E_{5+i,7})+
\sum_{r=1,3,4,5}\zeta_{5+r}(E_{2,5+r}-E_{r,7})],
\hspace{1.6cm}(6.57)\end{eqnarray*}
\begin{eqnarray*}T_3&=&-[\iota(E_{-\al_3}),T_2]=[(E_{3,2}-E_{7,8})|_{\cal
A}+(E_{3,2}-E_{7,8}),T_2]
\\&=&\zeta_3(D-c-6)+\sum_{i=1}^5\zeta_i(E_{3,i}-E_{5+i,8})
+ \sum_{r=1,2,4,5}\zeta_{5+r}(E_{3,5+r}-E_{r,8})],
\hspace{1.6cm}(6.58)\end{eqnarray*}
\begin{eqnarray*}T_4&=&-[\iota(E_{-\al_4}),T_3]=[(E_{4,3}-E_{8,9})|_{\cal
A}+(E_{4,3}-E_{8,9}),T_3]
\\&=&\zeta_4(D-c-6)+\sum_{i=1}^5\zeta_i(E_{4,i}-E_{5+i,9})
+ \sum_{r=1,2,3,5}\zeta_{5+r}(E_{4,5+r}-E_{r,9})],
\hspace{1.5cm}(6.59)\end{eqnarray*}
\begin{eqnarray*}T_5&=&-[\iota(E_{-\al_5}),T_4]=[(E_{5,4}-E_{9,10})|_{\cal
A}+(E_{5,4}-E_{9,10}),T_4]
\\&=&\zeta_5(D-c-6)+\sum_{i=1}^5\zeta_i(E_{5,i}-E_{5+i,10})
+\sum_{r=1}^4\zeta_{5+r}(E_{5,5+r}-E_{r,10}),\hspace{1.9cm}(6.60)\end{eqnarray*}
\begin{eqnarray*}T_{10}&=&-[\iota(E_{-\al_2}),T_4]=[(E_{10,4}-E_{9,5})|_{\cal
A}+(E_{10,4}-E_{9,5}),T_4]
\\&=&\zeta_{10}(D-c-6)+\sum_{i=1}^4\zeta_i(E_{10,i}-E_{5+i,5})
+\sum_{r=1}^5\zeta_{5+r}(E_{10,5+r}-E_{r,5}),\hspace{1.7cm}(6.61)\end{eqnarray*}
\begin{eqnarray*}T_9&=&-[\iota(E_{-\al_2}),T_5]=[(E_{10,4}-E_{9,5})|_{\cal
A}+(E_{10,4}-E_{9,5}),T_5]
\\&=&\zeta_9(D-c-6)+\sum_{i=1,2,3,5}\zeta_i(E_{9,i}-E_{5+i,4})
+\sum_{r=1}^5\zeta_{5+r}(E_{9,5+r}-E_{r,4}),\hspace{1.7cm}(6.62)\end{eqnarray*}
\begin{eqnarray*}T_8&=&[\iota(E_{-\al_4}),T_9]=[(E_{8,9}-E_{4,3})|_{\cal
A}+(E_{8,9}-E_{4,3}),T_9]
\\&=&\zeta_8(D-c-6)+\sum_{i=1,2,4,5}\zeta_i(E_{8,i}-E_{5+i,3})
+\sum_{r=1}^5\zeta_{5+r}(E_{8,5+r}-E_{r,3}),\hspace{1.7cm}(6.63)\end{eqnarray*}
\begin{eqnarray*}T_7&=&[\iota(E_{-\al_3}),T_8]=[(E_{7,8}-E_{3,2})|_{\cal
A}+(E_{7,8}-E_{3,2}),T_8]
\\&=&\zeta_7(D-c-6)+\sum_{i=1,3,4,5}\zeta_i(E_{7,i}-E_{5+i,2})
+\sum_{r=1}^5\zeta_{5+r}(E_{7,5+r}-E_{r,2}),\hspace{1.7cm}(6.64)\end{eqnarray*}
\begin{eqnarray*}T_6&=&[\iota(E_{-\al_1}),T_7]=[(E_{6,7}-E_{2,1})|_{\cal
A}+(E_{6,7}-E_{2,1}),T_7]
\\&=&\zeta_6(D-c-6)+\sum_{i=2}^5\zeta_i(E_{6,i}-E_{5+i,1})
+\sum_{r=1}^5\zeta_{5+r}(E_{6,5+r}-E_{r,1}).\hspace{2.3cm}(6.65)\end{eqnarray*}
Then ${\cal T}=\sum_{i=1}^{10}\mbb{C}T_i$ forms the 10-dimensional
natural module of $o(10,\mbb{C})$ with the standard basis
$\{T_1,...,T_{10}\}$.

Denote
$$T'_i=T_i-\zeta_i(D-c-6)\qquad\for\;\;i\in\ol{1,10}.\eqno(6.66)$$
Easily see that  ${\cal T}'=\sum_{i=1}^{10}\mbb{C}T_i'$ forms the
10-dimensional natural module of $o(10,\mbb{C})$ with the standard
basis $\{T_1',...,T_{10}'\}$. So we have the $o(10,\mbb{C})$-module
isomorphism from $U=\sum_{i=1}^{10}\mbb{C}\zeta_i$ to ${\cal T}'$
determined by $\zeta_i\mapsto T_i'$ for $i\in\ol{1,10}$. The weight
set of $U$ is
$$\Pi(U)=\{\pm\ves_1,...,\pm\ves_5\}.\eqno(6.67)$$
Let $\lmd\in\Lmd^+$. Denote
$$\Upsilon'(\lmd)=\{\lmd+\mu\mid
\mu\in\Pi(U),\;\;\lmd+\mu\in\Lmd^+\}.\eqno(6.68)$$ Take $M=V(\lmd)$.
It is known that
$$UV(\lmd)=U\otimes_{\mbb{C}}V(\lmd)\cong
\bigoplus_{\lmd'\in \Upsilon'(\lmd)}V(\lmd').\eqno(6.69)$$

Given $\lmd'\in\Upsilon'(\lmd)$, we pick a singular vector
$$u=\sum_{i=1}^{10}\zeta_iu_i\eqno(6.70)$$
of weight $\lmd'$ in $UV(\lmd)$, where $u_i\in V(\lmd)$. Moreover,
any singular vector of  weight $\lmd'$ in $UV(\lmd)$ is a scalar
multiple of $u$. Note that the vector
$$w=\sum_{i=1}^{10}T'_i(u_i)\eqno(6.71)$$
is also a  singular vector of  weight $\lmd'$ if it is not zero.
Thus
$$w=\flat_{\lmd'}u,\qquad \flat_{\lmd'}\in\mbb{C}.\eqno(6.72)$$
Set
$$\flat(\lmd)=\min\{\flat_{\lmd'}\mid\lmd'\in\Upsilon'(\lmd)\}.\eqno(6.73)$$
\psp

{\bf Theorem 6.4}. {\it The ${\cal G}^{E_6}$-module
$\widehat{V(\lmd)}$ is irreducible if}
$$c\in\mbb{C}\setminus\{\flat(\lmd)-6+\mbb{N},2\ell_\omega(\lmd)+2\mbb{N}\}.\eqno(6.74)$$

{\it Proof}. Recall that the ${\cal G}^{E_6}$-submodule $U({\cal
G}_-)(V(\lmd))$ is irreducible by Proposition 5.2. It is enough to
prove $\widehat{V(\lmd)}=U({\cal G}_-)(V(\lmd))$. It is obvious that
$$\widehat{V(\lmd)}_{\la 0\ra}=V(\lmd)=(U({\cal G}_-)(V(\lmd)))_{\la
0\ra}\eqno(6.75)$$ (cf. (6.3) and (6.4) with $M=V(\lmd)$). Moreover,
Lemma 6.1 with $M=V(\lmd)$,  (6.40) and (6.74) imply that
$\vf|_{\widehat{V(\lmd)}_{\la 1\ra}}$ is invertible, equivalently,
$$\widehat{V(\lmd)}_{\la 1\ra}=(U({\cal G}_-)(V(\lmd)))_{\la
1\ra}.\eqno(6.76)$$ Suppose that
$$\widehat{V(\lmd)}_{\la i\ra}=(U({\cal G}_-)(V(\lmd)))_{\la
i\ra}\eqno(6.77)$$ for $i\in\ol{0,k}$ with $1\leq k\in\mbb{N}$.

For any $v\in V(\lmd)$ and $\al\in\mbb{N}^{16}$ such that
$|\al|=k-1$, we have
$$T_r(x^\al v)=x^\al[(|\al|-c-6)\zeta_r+T'_r](v)\in (U({\cal G}_-)(V(\lmd)))_{\la
k+1\ra},\qquad r\in\ol{1,10}\eqno(6.78)$$ by (6.77) with $i=k-1,k$.
But
$$V'=\mbox{Span}\{[(|\al|-c-6)\zeta_r+T'_r](v)\mid r\in\ol{1,10},\;v\in
V(\lmd)\}\eqno(6.79)$$ forms an $o(10,\mbb{C})$-submodule of
$UV(\lmd)$ with respect to the action in (6.1). Let $u$ be a
$o(10,\mbb{C})$-singular in (6.70). Then
$$V'\ni
\sum_{r=1}^{10}[(|\al|-c-6)\zeta_r+T'_r](u_r)=(|\al|-c-6)u+w=(|\al|-c-6+\flat_{\lmd'})u\eqno(6.80)$$
by (6.71) and (6.72). Moreover, (6.73) and (6.74) yield $u\in V'$.
Since $UV(\lmd)$ is an $o(10,\mbb{C})$-module generated by all the
singular vectors, we have $V'=UV(\lmd)$. So
$$x^\al UV(\lmd)\subset (U({\cal G}_-)(V(\lmd)))_{\la
k+1\ra}.\eqno(6.81)$$ The arbitrariness of $\al$ implies
$$\zeta_r\widehat{V(\lmd)}_{\la k-1\ra}\subset (U({\cal G}_-)(V(\lmd)))_{\la
k+1\ra}\qquad\for\;\;r\in\ol{1,10}.\eqno(6.82)$$

Given any $f\in{\cal A}_k$ and $v\in V(\lmd)$, we have
$$\zeta_r\ptl_{x_i}(f)v\in \zeta_r\widehat{V(\lmd)}_{\la k-1\ra}\subset (U({\cal G}_-)(V(\lmd)))_{\la
k+1\ra}\qquad\for\;\;r\in\ol{1,10},\;i\in\ol{1,16}.\eqno(6.83)$$
Moreover,
\begin{eqnarray*}\qquad\eta_s(fv)&=&\iota(\eta_s)(fv)=P_s(fv)+f(\td\omega-c/2)(x_sv)\\
&\equiv&
f(k+\td\omega-c/2)(x_sv)\;\;(\mbox{mod}\;\sum_{r=1}^{10}\zeta_r\widehat{V(\lmd)}_{\la
k-1\ra})\hspace{4.1cm}(6.84)\end{eqnarray*} for $s\in\ol{1,16}$ by
(4.16)-(4.32), (5.17)-(5.32) and Lemma 6.1. According to (6.40),
(6.74), (6.82) and (6.84), we get
$$x_sfv\in (U({\cal G}_-)(V(\lmd)))_{\la
k+1\ra}\qquad\for\;\;s\in\ol{1,16}.\eqno(6.85)$$ Thus (6.77) holds
for $i=k+1$. By induction on $k$, (6.77) holds for any
$i\in\mbb{N}$, that is, $\widehat{V(\lmd)}=U({\cal
G}_-)(V(\lmd)).\qquad\Box$\psp

When $\lmd=0$, $V(0)$ is the one-dimensional trivial module and
$\ell_\omega(0)=\flat(0)=0$. So we have:\psp

{\bf Corollary 6.5}. {\it The ${\cal G}^{E_6}$-module
$\widehat{V(0)}$ is irreducible if}
$c\in\mbb{C}\setminus\{\mbb{N}-6\}.$ \psp

Next we consider the case $\lmd=k\ves_1=k\lmd_1$ for some positive
integer $k$, where $\lmd_1$ is the first fundamental weight. Note
$$\Upsilon(k\ves_1)=\{\lmd_4+k\ves_1,\lmd_4+(k-1)\ves_1+\ves_5\}\eqno(6.86)$$ by
(6.34) and (6.35). Thus (6.40) and (6.41) give
$$\ell_{\omega}(k\ves_1)=-4-k/2.\eqno(6.87)$$
In order to calculate $\flat(k\ves_1)$, we give a realization of
$V(k\ves_1)$. Observe that we have a representation of
$o(10,\mbb{C})$ on ${\cal B}=\mbb{C}[y_1,...,y_{10}]$ determined via
$$E_{i,j}|_{\cal
B}=y_i\ptl_{y_j}\qquad\for\;\;i,j\in\ol{1,10}.\eqno(6.88)$$ Denote
by ${\cal B}_k$ the subspace of homogenous polynomials in ${\cal B}$
with degree $k$. Set
$${\cal H}_k=\{h\in{\cal B}_k\mid
(\sum_{i=1}^5\ptl_{y_i}\ptl_{y_{5+i}})(h)=0\}.\eqno(6.89)$$ Then
${\cal H}_k\cong V(k\ves_1)$ and $y_1^k$ is a highest-weight vector.

According to (6.67) and (6.68),
$$\Upsilon'(k\ves_1)=\{(k+1)\ves_1,(k-1)\ves_1, k\ves_1+\ves_2\}.\eqno(6.90)$$
The vector $\zeta_1y_1^k$ is a singular vector in $U{\cal H}_k$ with
weight $(k+1)\ves_1$, where we take $M={\cal H}_k$ in the earlier
settings. By (6.55) and (6.66),
$$T_1'(y_1^k)=k\zeta_1y_1^k\lra \flat_{(k+1)\ves_1}=k.\eqno(6.91)$$
Moreover, $\zeta_1y_1^{k-1}y_2-\zeta_2y_1^k$  is a singular vector
in $U{\cal H}_k$ with weight $k\ves_1+\ves_2$. By (6.55), (6.57) and
(6.66), we find
\begin{eqnarray*}T'_1(y_1^{k-1}y_2)-T'_2(y_1^k)&=&(k-1)\zeta_1y_1^{k-1}y_2+\zeta_2y_1^k-k\zeta_1y_1^{k-1}y_2
\\ &=&\zeta_2y_1^k-\zeta_1y_1^{k-1}y_2=-(\zeta_1y_1^{k-1}y_2-\zeta_2y_1^k).\hspace{3.6cm}(6.92)\end{eqnarray*}
Thus $\flat_{k\ves_1+\ves_2}=-1$. Furthermore,
$$\varpi=(k+3)\sum_{i=1}^5[\zeta_iy_1^{k-1}y_{5+i}+\zeta_{5+i}y_1^{k-1}y_i]-(k-1)\zeta_1y_1^{k-2}\sum_{s=1}^5y_sy_{5+s}\eqno(6.93)$$
is a singular vector in $U{\cal H}_k$ with weight $(k-1)\ves_1$.
Expressions (6.55) and (6.57)-(6.66) yield
\begin{eqnarray*}& &(k+3)\sum_{i=1}^5[T'_i(y_1^{k-1}y_{5+i})+T'_{5+i}(
y_1^{k-1}y_i)]-(k-1)T_1'(x^\al y_1^{k-2}\sum_{s=1}^5y_sy_{5+s})
\\ &=&
(k+3)\sum_{i=1}^5[(k-1)\zeta_1
y_1^{k-2}y_iy_{5+i}-\sum_{r=1}^5\zeta_r y_1^{k-1}y_{5+r}
-\sum_{s\neq i}\zeta_{5+s} y_1^{k-1}y_s\\ & &
-((k-1)\dlt_{i,1}+1)(\sum_{r\neq i}\zeta_r
y_1^{k-1}y_{5+r}+\sum_{s=1}^5\zeta_{5+s} y_1^{k-1}y_s)\\
& &+(k-1)(1-\dlt_{i,1})\zeta_1y_1^{k-2}y_{5+i}y_i]
-(k-1)(k-2)\sum_{i=1}^5\zeta_1y^{k-2}_1y_iy_{5+i}
\\ &=&(-8-k)x^\al\varpi\lra
\flat_{(k-1)\ves_1}=-8-k.\hspace{7.4cm}(6.94)
\end{eqnarray*}
Therefore, $\flat(k\ves_1)=-8-k.$ By Theorem 6.4 and (6.87), we
obtain:\psp

 {\bf Corollary 6.6}. {\it The ${\cal G}^{E_6}$-module $\widehat{V(k\ves_1)}$ is
irreducible if} $c\in\mbb{C}\setminus\{\mbb{N}-14-k\}.$ \psp

Now we want to consider the cases $\lmd=\ves_1+\ves_2=\lmd_2$ (the
second fundamental weight) and $\lmd=\ves_1+\ves_2+\ves_3=\lmd_3$
(the third fundamental weight). Let ${\cal E}$ be the associative
algebra generated by $\{\sta_1,...,\sta_{10}\}$ with the defining
relations:
$$\sta_i\sta_j=-\sta_j\sta_i\qquad\for\;\;i,j\in\ol{1,10}.\eqno(6.95)$$
Moreover, we define linear operators
$\{\ptl_{\sta_1},...,\ptl_{\sta_{10}}\}$ on ${\cal E}$ by
$$\ptl_{\sta_i}(1)=0,\;\ptl_{\sta_i}(\sta_jw)=\dlt_{i,j}w-\sta_j\ptl_{\sta_i}(w)\qquad\for\;\;i,j\in\ol{1,10},\;
w\in{\cal E}.\eqno(6.96)$$ The representation of $o(10,\mbb{C})$ on
${\cal E}$ is defined via
$$E_{i,j}|_{\cal E}=\sta_i\ptl_{\sta_j}\qquad\for\;\;i,j\in\ol{1,10}.\eqno(6.97)$$
Denote
$${\cal E}_r=\mbox{Span}\{\sta_{i_1}\sta_{i_2}\cdots\sta_{i_r}\mid
1\leq i_1<i_2<\cdots<i_r\leq 10\}.\eqno(6.98)$$ Then ${\cal E}_2$
forms an irreducible $o(10,\mbb{C})$-submodule isomorphic
$V(\lmd_2)$ with a highest weight vector $\sta_1\sta_2$ and ${\cal
E}_3$ forms an irreducible $o(10,\mbb{C})$-submodule isomorphic
$V(\lmd_3)$ with a highest weight vector $\sta_1\sta_2\sta_3$.

Note that (6.35) gives
$$\Upsilon(\lmd_2)=\{\lmd_4+\lmd_2,\lmd_4,\lmd_4+\ves_1+\ves_5\},\eqno(6.99)$$
which yields $\ell_\omega(\lmd_2)=-8$. Moreover, (6.67) and (6.68)
imply
$$\Upsilon'(\lmd_2)=\{\ves_1+\lmd_2,\lmd_3,\ves_1\}.\eqno(6.100)$$
We have an $o(10,\mbb{C})$-singular vector $\zeta_1\sta_1\sta_2$ of
weight $\ves_1+\lmd_2$ in $U{\cal E}_2$, where we take $M={\cal
E}_2$ in the earlier settings. By (6.55) and (6.66),
$$T'_1(\sta_1\sta_2)=\zeta_1\sta_1\sta_2\lra\flat_{\ves_1+\lmd_2}=1\eqno(6.101)$$
Furthermore,,
$$u=\sum_{i=2}^5\zeta_{5+i}\sta_1\sta_i+\sum_{r=1}^5\zeta_r\sta_1\sta_{5+r}\eqno(6.102)$$
is an $o(10,\mbb{C})$-singular vector of weight $\ves_1$ in $U{\cal
E}_2$ and
$$w=\zeta_1\sta_2\sta_3-\zeta_2\sta_1\sta_3+\zeta_3\sta_1\sta_2\eqno(6.103)$$
is an $o(10,\mbb{C})$-singular vector of weight $\lmd_3$ in $U{\cal
E}_2$. According to (6.55), (6.57), (6.58) and (6.66), we have
\begin{eqnarray*}& &T'_1(\sta_2\sta_3)-T'_2(\sta_1\sta_3)+T'_3(\sta_1\sta_2)
\\&=&\zeta_2\sta_1\sta_3-\zeta_3\sta_1\sta_2-\zeta_1 \sta_2\sta_3-\zeta_3 \sta_1\sta_2-\zeta_1  \sta_2\sta_3+\zeta_2 \sta_1\sta_3
=-2w. \hspace{3.2cm}(6.104)\end{eqnarray*} So $\flat_{\lmd_3}=-2$.
Expressions (6.55) and  (6.57)-(6.66) give rise to
\begin{eqnarray*}& &\sum_{i=2}^5T'_{5+i}( \sta_1\sta_i)+\sum_{r=1}^5T'_r(
\sta_1\sta_{5+r})\\&=&-[\sum_{i=2}^5\zeta_{5+i}
\sta_1\sta_i+\sum_{r=2}^5\zeta_r
\sta_1\sta_{5+r}]-\sum_{i=2}^5[\sum_{i\neq
r\in\ol{1,5}}\zeta_r\sta_1\sta_{5+r}
+\sum_{ r=2}^5\zeta_{5+r}\sta_1\sta_r]\\
& & -\sum_{r=1}^5[\sum_{i=1}^5\zeta_i\sta_1\sta_{5+i}+\sum_{r\neq
i\in\ol{2,5}}\zeta_{5+i}\sta_1\sta_i]=-9u.\hspace{5.9cm}(6.105)
\end{eqnarray*}
Hence $\flat_{\ves_1}=-9$. Therefore, $\flat(\lmd_2)=-9$. Theorem
6.4 implies:\psp

 {\bf Corollary 6.7}. {\it The ${\cal G}^{E_6}$-module $\widehat{V(\lmd_2)}$ is
irreducible if} $c\in\mbb{C}\setminus\{\mbb{N}-16\}.$\psp

Observe that (6.35) gives
$$\Upsilon(\lmd_3)=\{\lmd_4+\lmd_3,\lmd_4+\lmd_2+\ves_5,
\lmd_4+\ves_1,\lmd_4+\ves_5\},\eqno(6.106)$$ which yields
$\ell_\omega(\lmd_3)=-21/2$. Moreover,
$$\Upsilon'(\lmd_2)=\{\ves_1+\lmd_3,\sum_{i=1}^4\ves_i,\lmd_2\}.\eqno(6.107)$$
Similarly we have $\flat_{\lmd_4+\lmd_3}=1$. Furthermore,,
$$u=\sum_{i=3}^5\zeta_{5+i}\sta_1\sta_2\sta_i+\sum_{r=1}^5\zeta_r\sta_1\sta_2\sta_{5+r}
\eqno(6.108)$$ is an $o(10,\mbb{C})$-singular vector of weight
$\lmd_2$ in $U{\cal E}_3$ and
$$w=\zeta_1\sta_2\sta_3\sta_4-\zeta_2\sta_1\sta_3\sta_4+\zeta_3\sta_1\sta_2\sta_4-\zeta_4\sta_1\sta_2\sta_3\eqno(6.109)$$
is an $o(10,\mbb{C})$-singular vector of weight $\sum_{i=1}^4\ves_i$
in $U{\cal E}_3$. According to (6.55), (6.57)-(6.59) and (6.66),
\begin{eqnarray*}&
&T'_1(\sta_2\sta_3\sta_4)-T_2'(\sta_1\sta_3\sta_4)+T'_3(\sta_1\sta_2\sta_4)-T_4'(\sta_1\sta_2\sta_3)
\\&=&\zeta_2\sta_1\sta_3\sta_4-\zeta_3\sta_1\sta_2\sta_4+\zeta_4\sta_1\sta_2\sta_3
-\zeta_1\sta_2\sta_3\sta_4-\zeta_3\sta_1\sta_2\sta_4+\zeta_4\sta_1\sta_2\sta_3
-\zeta_1\sta_2\sta_3\sta_4\\
& &+\zeta_2\sta_1\sta_3\sta_4+\zeta_4\sta_1\sta_2\sta_3
-\zeta_1\sta_2\sta_3\sta_4+\zeta_2\sta_1\sta_3\sta_4-\zeta_3\sta_1\sta_2\sta_4=-3w,
\hspace{2.6cm}(6.110)\end{eqnarray*}
 equivalently,
$\flat_{\sum_{i=1}^4\ves_i}=-3$. Expressions (6.55) and
(6.57)-(6.66) give rise to
\begin{eqnarray*}& &\sum_{i=3}^5T'_{5+i}(\sta_1\sta_2\sta_i)+\sum_{r=1}^5T_r(
\sta_1\sta_2\sta_{5+r})\\&=&-\sum_{i=3}^5(\zeta_{5+i}
\sta_1\sta_2\sta_i+\zeta_i\sta_1\sta_2\sta_{5+i})-\sum_{i=3}^5[\sum_{i\neq
r\in\ol{1,5}}\zeta_r\sta_1\sta_2\sta_{5+r}
+\sum_{ r=3}^5\zeta_{5+r}\sta_1\sta_2\sta_r]\\
& &
-\sum_{r=1}^5[\sum_{i=1}^5\zeta_i\sta_1\sta_2\sta_{5+i}+\sum_{r\neq
i\in\ol{3,5}}\zeta_{5+i}\sta_1\sta_2\sta_i]=-8u.\hspace{5.2cm}(6.111)
\end{eqnarray*}
So $\flat_{\lmd_2}=-8$. Therefore, $\flat(\lmd_3)=-8$. Theorem 6.4
yields:\psp

 {\bf Corollary 6.8}. {\it The ${\cal G}^{E_6}$-module $\widehat{V(\lmd_3)}$ is
irreducible if}
$c\in\mbb{C}\setminus\{\mbb{N}-15,-17,-19,-21\}.$\psp

Let $k$ be a positive integer. We calculate by (6.35) that
$$\Upsilon(k\lmd_4)=\{(k+1)\lmd_4,(k+1)\lmd_4-\ves_4+\ves_5,(k-1)\lmd_4+\ves_1\},\eqno(6.112)$$
which yields $\ell_\omega(k\lmd_4)=k/2-6$. The fifth fundamental
weight of $o(10,\mbb{C})$ is $\lmd_5=(1/2)\sum_{i=1}^5\ves_i$. Now
(6.35) implies
$$\Upsilon(k\lmd_5)=\{k\lmd_5+\lmd_4,k\lmd_5+\lmd_4-\ves_3-\ves_4,(k-1)\lmd_5\},\eqno(6.113)$$
which implies $\ell_\omega(k\lmd_5)=-k/2-10$.

We define a representation of $o(10,\mbb{C})$ on ${\cal
C}=\mbb{C}[z_1,...,z_{16}]$  obtained from (2.49)-(2.71) with ${\cal
A}$ replaced by ${\cal C}$ and $x_i$ replaced by $z_i$ for
$i\in\ol{1,16}$. Then the $o(10,\mbb{C})$-submodule ${\cal N}_k$
generated by $z_1^k$ is isomorphic to $V(k\lmd_4)$. Note that by
(6.67) and (6.68),
$$\Upsilon'(k\lmd_4)=\{k\lmd_4+\ves_1,k\lmd_4+\ves_5\}.\eqno(6.114)$$
Note that $\zeta_1z_1^k$ is an $o(10,\mbb{C})$-singular vector of
weight $k\lmd_4+\ves_1$ in $U{\cal N}_k$, where $M={\cal N}_k$ in
the earlier settings. By (2.71) with $x_i$ replaced by $z_i$, (6.55)
and (6.66),
$$T'_1(z_1^k)=\frac{k}{2}\zeta_1z_1^k\lra\flat_{k\lmd_4+\ves_1}=\frac{k}{2}.\eqno(6.115)$$
By (2.49)-(2.53), (2.69) and (2.70) with $x_i$ replaced by $z_i$,
$$(E_{9,10}-E_{5,4})(z_1^k)=kz_1^{k-1}z_2\in {\cal N}_k\lra z_1^{k-1}z_2\in {\cal
N}_k,\eqno(6.116)$$
$$(E_{8,9}-E_{4,3})(z_1^{k-1}z_2)=z_1^{k-1}z_3\in {\cal N}_k,\eqno(6.117)$$
$$(E_{3,2}-E_{7,8})(z_1^{k-1}z_3)=z_1^{k-1}z_4\in {\cal N}_k,\eqno(6.118)$$
$$(E_{2,1}-E_{6,7})(z_1^{k-1}z_4)=z_1^{k-1}z_7\in {\cal N}_k.\eqno(6.119)$$
According to (2.49)-(2.53)  with $x_i$ replaced by $z_i$ and Table
1, we find that the  vector
$$u'=\zeta_5z_1^k+\zeta_4z_1^{k-1}z_2+\zeta_3z_1^{k-1}z_3-\zeta_2z_1^{k-1}z_4+\zeta_1z_1^{k-1}z_7\eqno(6.120)$$
is an $o(10,\mbb{C})$-singular vector of weight $k\lmd_4+\ves_5$ in
$U{\cal N}_k$. Expressions (2.49)-(2.71) with $x_i$ replaced by
$z_i$, (6.55), (6.57)-(6.59) and (6.66),
\begin{eqnarray*}&
&T'_1( z_1^{k-1}z_7)-T_2'(z_1^{k-1}z_4)+T_3'( z_1^{k-1}z_3)+T'_4(
z_1^{k-1}z_2)+T_5'(z_1^k)
\\ &=&(k/2-1)\zeta_1 z_1^{k-1}z_7+\zeta_2
z_1^{k-1}z_4-\zeta_3 z_1^{k-1}z_3-\zeta_4 z_1^{k-1}z_2-\zeta_5
z_1^k\\ & &-\zeta_1z_1^{k-1}z_7-(k/2-1)\zeta_2 z_1^{k-1}z_4 -\zeta_3
z_1^{k-1}z_3-\zeta_4 z_1^{k-1}z_2-\zeta_5 z_1^k
\\ & &-\zeta_1z_1^{k-1}z_7+\zeta_2
z_1^{k-1}z_4+(k/2-1)\zeta_3 z_1^{k-1}z_3-\zeta_4
z_1^{k-1}z_2-\zeta_5 z_1^k\\ & &-\zeta_1z_1^{k-1}z_7+\zeta_2
z_1^{k-1}z_4-\zeta_3 z_1^{k-1}z_3+(k/2-1)\zeta_4
z_1^{k-1}z_2-\zeta_5 z_1^k\\ & &-k\zeta_1z_1^{k-1}z_7+k\zeta_2
z_1^{k-1}z_4-k\zeta_3 z_1^{k-1}z_3-k\zeta_4
z_1^{k-1}z_2-(k/2)\zeta_5 z_1^k\\
&=&-(k/2+4)u',\hspace{11.2cm}(6.121)
\end{eqnarray*}
equivalently, $\flat_{k\lmd_4+\ves_5}=-(k/2+4)$. Thus
$\flat(k\lmd_4)=-(k/2+4)$. Symmetrically, $\flat(k\lmd_5)=-(k/2+4)$.
By Theorem 6.4, we have:\psp

{\bf Corollary 6.9}. {\it The ${\cal G}^{E_6}$-module
$\widehat{V(k\lmd_4)}$ is irreducible if}
$c\in\mbb{C}\setminus\{\mbb{N}-10-k/2,2\mbb{N}+k-12\}.$ {\it The
${\cal G}^{E_6}$-module $\widehat{V(k\lmd_5)}$ is irreducible if}
$c\in\mbb{C}\setminus\{\mbb{N}-10-k/2,2\mbb{N}-k-20\}$.

\vspace{1cm}

\noindent{\Large \bf References}

\hspace{0.5cm}

\begin{description}

\item[{[A]}] J. Adams, {\it Lectures on Exceptional Lie Groups}, The
University of Chicago Press Ltd., London, 1996.

\item[{[AB1]}] G. Anderson and T. Bla\v{z}ek, $E_6$ unification model
building.I. Clebsch-Gordan coefficients of $27\otimes\ol{27}$, {\it
J. Math. Phys.} {\bf 41} (2000), no. 7, 4808-4816.

\item[{[AB2]}] G. Anderson and T. Bla\v{z}ek, $E_6$ unification model
building.II. Clebsch-Gordan coefficients of $78\otimes 78$, {\it J.
Math. Phys.} {\bf 41} (2000), no. 12, 8170-8189.

\item[{[AB3]}] G. Anderson and  T. Bla\v{z}ek, $E_6$ unification model
building.III. Clebsch-Gordan coefficients in $E_6$ tensor products
of the 27 with higher-dimensional representations, {\it J. Math.
Phys.} {\bf 46} (2005), no. 1, 013506, 13pp.

\item[{[As1]}] M. Aschbacher, The 27-dimensional module for
$E_6$.I., {\it Invent. Math.} {\bf 89} (1987), no. 1, 159-195.

\item[{[BCDH]}] P. Berglund, P. Candelas, X. de le Ossa, E. Derrick,
J. Distler and T. H\"{u}bsch, On instanton contrbutions to the
masses and couplings of $E_6$ singles, {\it Nuclear Phys. B} {\bf
454} (1995), no. 1-2, 127-163.

\item[{[BZ]}] B. Bineger and R. Zierau, A singular representation of
$E_6$, {\it Trans. Amer. Math. Soc.} {\bf 341} (1994), no. 2,
771-785.

\item[{[B-N]}] J. Bion-Nadal, Subfactor of the hyperfinite $\Pi_1$
factor with Coxeter graph $E_6$ as invariant, {\it J. Operator
Theory} {\bf 28} (1992), 27-50.

\item[{[BK]}] R. Brylinski and B. Kostant, Minimal representations
of $E_6,\;E_7$, and $E_8$ and the generalized Capelli identity, {\it
Proc. Nar. Acad. Sci. U.S.A.} {\bf 91} (1994), no. 7., 2469-2472.

\item[{[CD]}] B. Cerchiai and A. Scotti, Mapping the geometry of the
$E_6$ group, {\it J. Math. Phys.} {\bf 49} (2008), no. 1, 012107,
19pp.

\item[{[DL]}] C. Das and L. Laperashvili, Preon model and family
replicated $E_6$ unification, {\it SIGMA} {\bf 4} (2008), 012,
15pages.

\item[{[D]}] L. Dickson, A class of groups in an arbitrary realm
connected with the configuration of the 27 lines on a cubic surface,
{\it J. Math.} {\bf 33} (1901), 145-123.

\item[{[FGP]}] J. Fern\'{a}ndez-N\'{u}$\td{a}$ez, W. Garcia-Fuertes
and A. Perelomov,  Irreducible characters and Clebsch-Gordan series
for the exceptional algebra $E_6$: an approach through the quantum
Calogero-Sutherland model, {\it J. Nonlinear Math. Phys.} {\bf 12},
suppl. 1, 280-301.

\item[{[Gm]}] M. D. Gould, Tensor operators and projection techniques in infinite dimensional
  representations of semi-simple Lie algebras, {\it J. Phys. A: Math. Gen.} {\bf 17} (1984), 1-17.

\item[{[GSA]}] A. Ghezelbash, A. Shafiekhani and M. Abolbasani, On
the Picard-Fuchs equations of $N=2$ supersymmetric $E_6$ Yang-Mills
theory, {\it Modern Phys. Lett. A} {\bf 13} (1998), no. 7, 527-531.

\item[{[G]}] D. Ginzburg, On standard $L$-functions for $E_6$ and
$E_7$, {\it J. Reine Angew. Math.} {\bf 465} (1995), 101-131.

\item[{[HM]}] N. Haba and T. Matsuoka, Large lepton flavor mixing
and $E_6$-type unification models, {\it Progr. Theoret. Phys.} {\bf
99} (1998), no. 5, 831-842.

\item[{[HH]}] J. E. M. Homos and Y. M. M. Homos, Algebraic model for
the evolution of the generic code, {\it Phys. Rev. Lett.} {\bf 71}
(1991), 4401-4404.

\item[{[Hub]}] R. Hubert, The $(A_2, G_2)$ duality in $E_6$,
octonions and the triality principle, {\it Trans. Amer. Math. Soc.}
{\bf 360} (2008), no. 1, 347-367.

\item[{[Hum]}] J. E. Humphreys, {\it Introduction to Lie Algebras and Representation Theory},
 Springer-Verlag New York Inc., 1972.

\item[{[HK]}] R. Howl and S. King, Minimal $E_6$ supersymmetric
standard model, {\it J. High Ener. Phys.} {\bf 01}(2008), 030, 31pp.

\item[{[I]}] A. Iltyakov, On rational invariants of the group $E_6$,
{\it Proc. Math. Soc.} {\bf 124} (1996), no. 12, 3637-3640.

\item[{[Ka]}] V. Kac, {\it Infinite-Dimensional Lie Algebras},
Birkh\"{a}ser, Boston, Inc., 1982.

\item[{[K]}] B. Kostant, On the tensor product of a finite and an infinite
dimensional representation, {\it J. Func. Anal.} {\bf 20} (1975),
257-285.

\item[{[MPW]}] I. Morrison, P. Pieruschka  and
B. Wybourne, The interacting boson model with the exceptional groups
$G_2$ and $E_6$, {\it J. Math. Phys.} {\bf 32} (1991), no. 2,
356-372.

\item[{[OM]}] Y. Okamoto and R. Marshak, A garnd unification preon
model with $E_6$ metacolor, {\it Nuclear Phys. B} {\bf 268} (1986),
no. 2, 397-405.

\item[{[R]}] H. Rubenthaler, The $(A_2,G_2)$ duality in $E_6$,
octonians and the triality principle, {\it Trans. Amer. Math. Soc.}
{\bf 360} (2008), no. 1, 347-367.

\item[{[S1]}] G. Shen, Graded modules of graded Lie algebras of
Cartan type (I)---mixed product of modules, {\it Science in China A}
{\bf  29} (1986), 570-581.

\item[{[S2]}] G. Shen, Graded modules of graded Lie algebras of
Cartan type (II)---positive and negative graded modules, {\it
Science in China A} {\bf  29} (1986), 1009-1019.

\item[{[S3]}] G. Shen, Graded modules of graded Lie algebras of
Cartan type (III)---irreducible modules, {\it Chin. Ann. of Math B}
{\bf  9} (1988), 404-417.

\item[{[SW]}] K. Suzuki and M. Wakui, On the Turaev-Viro-Ocneanu
invariant of 3-manifolds derived from the $E_6$-subfactor, {\it
Kyushu J. Math}. {\bf 56} (2002), 59-81.

\item[{[T]}] J. Tits, A local approach to buildings, in: C. Davis, B. Gronbaum  and F. Sherk (eds),
   ``Geometric Vein,"  Berlin-Heidelberg-New York, Springer, 1982, pp.
519-547.

\item[{[Wa]}] X. Wang, Identification of Gepner's model with twisted
LG model and $E_6$ singlets, {\it Modern Phy. Lett. A} {\bf 6}
(1991), no. 23, 2155-2162.

\item[{[X1]}] X. Xu, {\it Kac-Moody Algebras and Their Representations}, China
Science Press, 2007.

\item[{[X2]}] X. Xu, A cubic $E_6$-generalization of the classical theorem on harmonic
polynomials, {\it J. Lie Theory} {\bf 21} (2011), 145-164.

\item[{[X3]}] X. Xu, Representations of Lie algebras and coding
theory, {\it J. Lie Theory}, accepted.

\item[{[XZ]}] X. Xu and Y. Zhao, Generalized conformal representations
of orthogonal Lie algebras, {\it arXiv:1105.1254v1[math.RT].}

\item[{[ZX]}] Y. Zhao and X. Xu, Generalized  projective  representations
 for sl(n+1), {\it J. Algebra} {\bf 328} (2011), 132-154.

\end{description}

\end{document}